\documentclass[12pt,reqno]{amsart}
\usepackage{amssymb}
\usepackage{amsxtra}
\usepackage{mathrsfs}
\usepackage[all]{xy}
\newtheorem{theorem}{Theorem}[section]
\newtheorem*{theorem*}{Theorem}
\newtheorem{lemma}[theorem]{Lemma}
\newtheorem{prop}[theorem]{Proposition}
\newtheorem{corollary}[theorem]{Corollary}

\theoremstyle{definition}
\newtheorem{definition}{Definition}[section]
\theoremstyle{remark}
\newtheorem{example}{Example}[section]
\newtheorem{remark}{Remark}[section]
\renewcommand*{\Im}{\mathop{\mathrm{Im}}}

\newcommand*{\h}{\mathbf h}

\DeclareMathOperator{\Tor}{Tor}

\DeclareMathOperator{\Ker}{Ker}
\newcommand*{\Bil}{\mathscr Bil}
\DeclareMathOperator{\Coker}{Coker}
\DeclareMathOperator{\coker}{coker}
\newcommand*{\ptens}[1]{\mathop{\widehat\otimes}_{#1}}

\newcommand*{\Ptens}{\mathop{\widehat\otimes}}

\newcommand*{\id}{1}
\newcommand*{\op}{\mathrm{op}}

\DeclareMathOperator{\rann}{rann}
\newcommand*{\un}{\mathrm{un}}
\newcommand*{\ann}{\mathrm{ann}}
\newcommand*{\lmod}{\mbox{-}\!\mathop{\mathsf{mod}}}
\newcommand*{\rmod}{\mathop{\mathsf{mod}}\!\mbox{-}}
\newcommand*{\bimod}{\mbox{-}\!\mathop{\mathsf{mod}}\!\mbox{-}}
\newcommand*{\lunmod}{\mbox{-}\!\mathop{\mathsf{unmod}}}
\newcommand*{\runmod}{\mathop{\mathsf{unmod}}\!\mbox{-}}
\newcommand*{\biunmod}{\mbox{-}\!\mathop{\mathsf{unmod}}\!\mbox{-}}
\newcommand*{\lbarmod}{\mbox{-}\!\underline{\mathop{\mathsf{mod}}}}
\newcommand*{\rbarmod}{\underline{\mathop{\mathsf{mod}}}\!\mbox{-}}
\newcommand*{\bibarmod}{\mbox{-}\!\underline{\mathop{\mathsf{mod}}}\!\mbox{-}}
\newcommand*{\lunbarmod}{\mbox{-}\!\underline{\mathop{\mathsf{unmod}}}}
\newcommand*{\runbarmod}{\underline{\mathop{\mathsf{unmod}}}\!\mbox{-}}
\newcommand*{\biunbarmod}{\mbox{-}\!\underline{\mathop{\mathsf{unmod}}}\!\mbox{-}}
\newcommand*{\lcmod}{\mbox{-}\!\mathop{\mathsf{lcmod}}}
\newcommand*{\LCS}{\mathsf{LCS}}
\newcommand*{\hLCS}{\mathsf{LCS}\sphat\; }
\newcommand*{\Ban}{\mathsf{Ban}}
\newcommand*{\Fr}{\mathsf{Fr}}
\newcommand*{\Vect}{\mathsf{Vect}}
\newcommand*{\CC}{\mathbb C}
\newcommand*{\cC}{\mathscr C}
\newcommand*{\cG}{\mathscr G}
\newcommand*{\cX}{\mathscr X}
\newcommand*{\cL}{\mathscr L}
\newcommand*{\cH}{\mathscr H}
\newcommand*{\A}{\mathscr A}
\newcommand*{\B}{\mathscr B}
\newcommand*{\cE}{\mathscr E}

\newcommand*{\N}{\mathbb N}
\newcommand*{\Z}{\mathbb Z}

\newcommand*{\lar}{\leftarrow}

\newcommand*{\xla}{\xleftarrow}
\newcommand*{\xra}{\xrightarrow}

\newcommand*{\la}{\langle}
\newcommand*{\ra}{\rangle}
\newcommand*{\eps}{\varepsilon}

\begin{document}
\title[Flat cyclic Fr\'echet modules]%
{Flat cyclic Fr\'echet modules,\\ amenable Fr\'echet algebras,
\\ and approximate identities}
\author{A. Yu. Pirkovskii}
\address{Department of Nonlinear Analysis and Optimization\\
Faculty of Science\\
Peoples' Friendship University of Russia\\
Mikluho-Maklaya 6\\
117198 Moscow\\
Russia}
\email{pirkosha@sci.pfu.edu.ru, pirkosha@online.ru}
\thanks{Partially supported by the RFBR grants 05-01-00982 and 05-01-00001,
and by the President of Russia grant MK-1490.2006.1.}
\begin{abstract}
Let $A$ be a locally $m$-convex
Fr\'echet algebra. We give a necessary and sufficient condition
for a cyclic Fr\'echet $A$-module $X=A_+/I$ to be strictly flat,
generalizing thereby a criterion of Helemskii and Sheinberg.
To this end, we introduce a notion of ``locally bounded approximate identity''
(a locally b.a.i. for short),
and we show that $X$ is strictly flat if and only if the ideal $I$
has a right locally b.a.i. Next we apply this result to amenable
algebras and show that a locally $m$-convex
Fr\'echet algebra $A$ is amenable if and only if $A$ is
isomorphic to a reduced inverse limit of amenable Banach algebras.
We also extend a number of characterizations of amenability obtained by Johnson and by
Helemskii and Sheinberg to the setting of locally $m$-convex Fr\'echet algebras.
As a corollary, we show that Connes and Haagerup's theorem on
amenable $C^*$-algebras and Sheinberg's theorem on amenable uniform
algebras hold in the Fr\'echet algebra case.
We also show that a quasinormable locally $m$-convex Fr\'echet algebra
has a locally b.a.i. if and only if it has a b.a.i.
On the other hand, we give an example of a commutative, locally
$m$-convex Fr\'echet-Montel algebra which has a locally b.a.i.,
but does not have a b.a.i.
\end{abstract}
\maketitle

\section{Introduction}

The notion of ``flat Banach module'' over a Banach algebra
was introduced by A.~Ya.~Helemskii \cite{Hel_period} within the framework
of his ``Topological Homology'' theory, which is a version
of relative homological algebra in categories of topological modules
over topological algebras (see \cite{X1,X_HOA}).
By definition, a left Banach module $X$ over a Banach
algebra $A$ is flat if the projective tensor product
$(\,\cdot\,)\ptens{A} X$ takes ``admissible'' exact sequences of right
Banach $A$-modules to exact sequences of vector spaces (for details,
see Section~\ref{sect:flat}). About the same time B.~E.~Johnson
published his important memoir \cite{Jhnsn_CBA} where he introduced amenable
Banach algebras. By definition, a Banach algebra $A$ is amenable if
for each Banach $A$-bimodule $X$ every continuous derivation from $A$
to the dual bimodule, $X^*$, is inner.
A relation between flat modules and amenable algebras was discovered
by A.~Ya.~Helemskii and M.~V.~Sheinberg \cite{Hel_Shein}.
They proved that a Banach algebra $A$ is amenable if and only if $A_+$,
the unitization of $A$, is a flat Banach $A$-bimodule.
Using this, they gave several useful characterizations of amenability
in terms of bounded approximate identities (b.a.i.'s); see Theorem~\ref{thm:Hel_amen}
below. These characterizations were obtained as a consequence of
the following flatness criterion for cyclic Banach modules
\cite{Hel_period, Hel_Shein}
(for terminology, see Section~\ref{sect:flat}).

\begin{theorem}[Helemskii, Sheinberg]
\label{thm:Hel_flat}
Let $A$ be a Banach algebra, and let $I\subset A_+$ be a closed left
ideal. Then the following conditions are equivalent:

\begin{itemize}
\item[{\upshape (i)}] $A_+/I$ is strictly flat;
\item[{\upshape (ii)}] $I$ has a right b.a.i.
\end{itemize}

If, in addition, $I$ is weakly complemented in $A_+$
(i.e., if the annihilator of $I$ is complemented in $A_+^*$),
then {\upshape (i)} and {\upshape (ii)} are equivalent to

\begin{itemize}
\item[{\upshape (iii)}] $A_+/I$ is flat.
\end{itemize}
\end{theorem}

In order to formulate the above-mentioned characterizations of
amenability obtained by Helemskii and Sheinberg, let us introduce
some notation. Given a Banach algebra $A$, set $A^e=A_+\Ptens A_+^\op$
(respectively, $A^e_0=A\Ptens A^\op$), where $A^\op$ is the algebra opposite to $A$,
and let $I^\Delta$ (respectively, $I_0^\Delta$) denote the kernel
of the product map $A^e\to A_+$ (respectively, $A^e_0\to A$).

\begin{theorem}[Helemskii, Sheinberg]
\label{thm:Hel_amen}
Let $A$ be a Banach algebra. The following conditions are equivalent:
\begin{enumerate}
\renewcommand{\theenumi}{\roman{enumi}}
\item $A$ is amenable;
\item $A_+$ is a strictly flat Banach $A$-bimodule;
\item $A$ is biflat and has a b.a.i.;
\item $I^\Delta$ has a right b.a.i.;
\item $A$ has a b.a.i., and $I^\Delta_0$ has a right b.a.i.
\end{enumerate}
\end{theorem}

Condition (v) of Theorem~\ref{thm:Hel_amen} is similar in spirit
to an earlier result of Johnson \cite{Jhnsn_appr}
(for terminology, see Section~\ref{sect:flat_amen}).

\begin{theorem}[Johnson]
\label{thm:Jhnsn}
Let $A$ be a Banach algebra. The following conditions are equivalent:
\begin{enumerate}
\renewcommand{\theenumi}{\roman{enumi}}
\item $A$ is amenable;
\item $A$ has a bounded approximate diagonal;
\item $A$ has a virtual diagonal.
\end{enumerate}
\end{theorem}

The aim of the present paper is to generalize Theorems~\ref{thm:Hel_flat},
\ref{thm:Hel_amen}, and~\ref{thm:Jhnsn} to locally $m$-convex Fr\'echet algebras
(i.e., to Fr\'echet-Arens-Michael algebras, in the terminology of \cite{X2}).
Let us remark that the question of whether Theorem~\ref{thm:Hel_amen}
holds for nonnormable locally convex algebras was posed explicitly in \cite{X_31}.
The basic difficulty here is that the duality between flat and injective
Banach modules, which was heavily used by Helemskii and Sheinberg,
is not available within the framework of Fr\'echet modules.
In fact, given a Fr\'echet module over a Fr\'echet algebra $A$, the dual
module, $X^*$, has no reasonable topology making it into a Fr\'echet space.
Moreover, the action of $A$ on $X^*$ often fails to be jointly continuous
with respect to any ``natural'' topology on $X^*$ (cf. \cite{T1}).
Let us also remark that many Fr\'echet algebras do not have nontrivial
injective Fr\'echet modules at all \cite{Pir_inject, Pir_injdim, Pir_nova}.

To overcome this difficulty, we use the Arens-Michael decomposition
for complete locally $m$-convex algebras \cite{Arens_gennorm,Michael} and for
complete locally convex modules over such algebras \cite{Pir_qfree}.
It turns out, however, that the property of having a bounded approximate identity
does not behave well under the Arens-Michael decomposition.
That is why we have to introduce a new notion of ``locally bounded approximate
identity'', which might be of independent interest.

The paper is organized as follows. Section~\ref{sect:prelim} contains
some background material on topological algebras and modules.
In Section~\ref{sect:tensprod} we prove some
general results on the behavior of the projective tensor product
functor. In particular, we show that this functor is cokernel-preserving.
In Section~\ref{sect:flat} we discuss flat and strictly flat Fr\'echet modules,
and we show that a (strictly) flat Banach module over a Banach algebra remains (strictly)
flat when considered as a Fr\'echet module. In Section~\ref{sect:bai}
we prove some general facts on approximate identities and approximate diagonals
in topological algebras.
Although these facts are well known in the case of Banach
algebras, we give full proofs for the sake of completeness.
In Section~\ref{sect:lbai} we introduce locally bounded approximate identities
and locally bounded approximate diagonals and study their basic properties.
Section~\ref{sect:strflat} contains the main result of the paper
(Theorem~\ref{thm:strflat} and Corollary~\ref{cor:strflat}), which extends
the first part of Helemskii and Sheinberg's Theorem~\ref{thm:Hel_flat}
(specifically, the equivalence of (i) and (ii))
to Fr\'echet-Arens-Michael algebras.
In Section~\ref{sect:qsinorm}, by using a result of Palamodov~\cite{Pal_methods},
we show that for quasinormable Fr\'echet-Arens-Michael algebras
the notions of ``locally bounded'' and ``bounded'' approximate identities
are equivalent. In Section~\ref{sect:flat_amen} we extend the second part
of Helemskii and Sheinberg's Theorem~\ref{thm:Hel_flat}
(specifically, the equivalence of (i) and (iii)) to
Fr\'echet-Arens-Michael algebras, and then apply this result to characterizing
amenability for such algebras, in the spirit of Theorems~\ref{thm:Hel_amen}
and~\ref{thm:Jhnsn}. In particular, we obtain a partial answer to a question
posed by Helemskii \cite[Problem 11]{X_31}.
We also prove that a Fr\'echet-Arens-Michael algebra
is amenable if and only if it is isomorphic to a reduced inverse limit of
amenable Banach algebras. This result is then used to show
that a $\sigma$-$C^*$-algebra is amenable
if and only if it is nuclear, and that a uniform Fr\'echet algebra is amenable
if and only if it is isomorphic to the algebra of continuous
functions on a hemicompact $k$-space.
Finally, in Section~\ref{sect:counterexample}
we give an example of a (non-quasinormable)
commutative Fr\'echet-Arens-Michael algebra (which is, in addition, a Montel space)
with a locally b.a.i., but without a b.a.i.
This example shows that Helemskii and Sheinberg's
Theorem~\ref{thm:Hel_flat} does not extend {\em verbatim} to
locally $m$-convex Fr\'echet algebras.

Some of the results of this paper were announced in \cite{Pir_ICTAA05}.
A related result was independently obtained by C.~P.~Podara \cite{Podara}.
Specifically, she proves that implication $\mathrm{(ii)}\Rightarrow\mathrm{(i)}$
of Helemskii and Sheinberg's Theorem~\ref{thm:Hel_flat} holds for every
(not necessarily locally $m$-convex) Fr\'echet algebra, and that
condition (ii) implies also that $I$ is strictly flat.

\section{Preliminaries}
\label{sect:prelim}
We shall work over the complex numbers $\CC$. Let $\Vect$ denote the
category of all vector spaces and linear maps.
The category
of locally convex spaces (l.c.s.'s) and continuous linear maps
will be denoted by $\LCS$.
Let us remark that we do not require locally convex spaces
to be Hausdorff, but we follow the convention that complete l.c.s.'s
are Hausdorff by definition.
The full subcategories of $\LCS$ consisting of complete
l.c.s.'s (respectively, of Fr\'echet spaces, of Banach spaces)
will be denoted by $\hLCS$ (respectively, $\Fr$, $\Ban$).

The completion of an l.c.s. $E$ is defined to be the completion of
the associated Hausdorff l.c.s., $E/\overline{\{ 0\}}$, and is denoted by
$E\sptilde$. We denote by $E^*$ the dual of $E$ and always endow $E^*$
with the strong topology, unless stated otherwise.
The canonical embedding of $E$ into $E^{**}$ is denoted by $i_E$.
The space of all continuous linear maps between l.c.s.'s
$E$ and $F$ is denoted by $\cL(E,F)$. Given l.c.s.'s $E,F,G$, we denote by
$\Bil(E\times F,G)$ the space of all jointly continuous
bilinear maps from $E\times F$ to $G$.
The completed projective tensor product of $E$ and $F$ is denoted by $E\Ptens F$.

By a {\em topological algebra} we mean a topological vector space $A$
together with the structure of an associative algebra such that the product
map $A\times A\to A$ is separately continuous. In what follows, when
using the word ``algebra'' with an adjective that describes a linear topological
property (such as ``locally convex'', ``complete'', ``Fr\'echet'' etc.), we mean that
the underlying topological vector space of the algebra in question has the
specified property. The same applies to topological modules (see below).
A complete locally convex algebra
with jointly continuous product is called a {\em $\Ptens$-algebra} \cite{T1,X1}.
A seminorm $\|\cdot\|$ on an algebra $A$ is {\em submultiplicative} if
$\| ab\|\le \| a\| \| b\|$ for all $a,b\in A$.
A locally convex algebra $A$ is {\em locally $m$-convex}
\cite{Michael} if the topology on $A$ can be defined by a family
of submultiplicative seminorms.
Equivalently, $A$ is locally $m$-convex if
it has a base $\mathscr U$ of $0$-neighborhoods consisting of absolutely convex,
idempotent sets (here ``idempotent'' means that $U^2\subset U$ for each $U\in\mathscr U$).
An {\em Arens-Michael algebra} \cite{X2} is a complete locally $m$-convex
algebra. A {\em Fr\'echet-Arens-Michael algebra} is a locally $m$-convex
Fr\'echet algebra.

Let $A$ be a topological algebra. Following \cite{Jhnsn_centr}, we say that
$A$ is {\em left hypotopological} if the product $A\times A\to A$ is left hypocontinuous
w.r.t. the family of all bounded subsets of $A$.
This means that for each $0$-neighborhood $V\subset A$ and each bounded subset
$B\subset A$ there exists a $0$-neighborhood $U\subset A$ such that $UB\subset V$.
Right hypotopological algebras are defined similarly.
``Hypotopological'' means ``left and right hypotopological''.
It is easy to show (cf. \cite{Jhnsn_centr}) that
the product of two bounded sets in a left hypotopological algebra is bounded.
Indeed, given bounded sets $B_1,B_2\subset A$ and a $0$-neighborhood
$U\subset A$, we may successively find a $0$-neighborhood $V\subset A$ such that
$VB_2\subset U$ and $\lambda>0$ such that $B_1\subset\lambda V$.
Then $B_1 B_2\subset\lambda U$, so that $B_1 B_2$ is bounded.
The same argument works in the case of a right hypotopological algebra.

If $A$ is a left hypotopological locally convex algebra, then the strong bidual,
$A^{**}$, can be endowed with a product making it into a topological
algebra \cite{Gulick}. The details of the construction are similar to the
Banach algebra case (see, e.g., \cite[2.6]{Dales}, \cite[1.4]{Palmer}). Specifically,
given $f\in A^*$ and $a\in A$, define $f\cdot a\in A^*$ by
$\la f\cdot a,b\ra=\la f,ab\ra\; (b\in A)$. Next, given $u\in A^{**}$
and $f\in A^*$, define $u\cdot f\in A^*$ by
$\la u\cdot f,a\ra=\la u,f\cdot a\ra\; (a\in A)$.
Finally, given $u,v\in A^{**}$, define $uv\in A^{**}$ by
$\la uv,f\ra=\la u,v\cdot f\ra\; (f\in A^*)$.
The bilinear map $(u,v)\mapsto uv$ is a separately continuous, associative
product on $A^{**}$ ({\em the first Arens product}) \cite[3.8, 3.9]{Gulick}.
The canonical embedding $i_A\colon A\to A^{**}$ is an algebra homomorphism
\cite[3.5]{Gulick}.
The first Arens product is weak$^*$-weak$^*$ continuous w.r.t. the first variable,
i.e., for each $v\in A^{**}$ the map $A^{**}\to A^{**},\; u\mapsto uv$, is continuous
w.r.t. the weak$^*$ topology on both copies of $A^{**}$ \cite[3.4]{Gulick}.
In general, the first Arens product does not have the above property w.r.t. the
second variable. However, for each $a\in A$ the map
$A^{**}\to A^{**},\; v\mapsto i_A(a)v$ is continuous
w.r.t. the weak$^*$ topology on both copies of $A^{**}$ \cite[3.6]{Gulick}.

Given a topological algebra $A$, a {\em left semitopological $A$-module}
is a topological vector space $X$ together with the structure of a left
$A$-module such that for each $a\in A$ the map
$X\to X,\; x\mapsto a\cdot x$ is continuous.
If, in addition, for each $x\in X$ the map $A\to X,\; a\mapsto a\cdot x$ is continuous,
then we say that $X$ is a {\em left topological $A$-module}.
Right (semi)topological
$A$-modules and (semi)topological $A$-bimodules are defined similarly.
In what follows, ``locally convex $A$-module'' means
``locally convex {\em topological} $A$-module''.
If $A$ is a $\Ptens$-algebra and $X$ is a left locally
convex $A$-module, we say that $X$ is a left {\em $A$-$\Ptens$-module} \cite{T1,X1}
if $X$ is complete and the action $A\times X\to X$ is jointly continuous.
If $X$ and $Y$ are left semitopological $A$-modules
(respectively, right semitopological $A$-modules, semitopological $A$-bimodules),
then the space of all continuous left $A$-module morphisms
(respectively, right $A$-module morphisms, $A$-bimodule morphisms)
from $X$ to $Y$ is denoted by ${_A}\h(X,Y)$
(respectively, $\h_A(X,Y)$, ${_A}\h_A(X,Y)$).

If $A$ is a locally convex algebra and $X$ is a left semitopological locally convex
$A$-module, then $X^*$ is a right semitopological locally convex $A$-module w.r.t. the
action $\la f\cdot a,x\ra=\la f,a\cdot x\ra\; (f\in X^*,\; a\in A,\; x\in X)$.
Similarly, if $X$ is a right semitopological locally convex $A$-module, then $X^*$ is
a left semitopological locally convex $A$-module.
If the action of $A$ on $X$ is hypocontinuous w.r.t.
the family of all bounded subsets of $X$,
then the action of $A$ on $X^*$ is separately continuous, so that
$X^*$ is a locally convex $A$-module (cf. \cite[Section 3]{T1} or \cite[3.1]{Gulick}).
Note that the above condition is satisfied automatically
whenever $A$ is barrelled and $X$ is a locally convex $A$-module.

Let $A$ be an algebra, and let $A_+$ denote the unitization of $A$.
We denote by $j_A$ the embedding of $A$ into $A_+$, and by
$\eps_A\colon A_+\to\CC$ the homomorphism that vanishes on $A$
and satisfies $\eps_A(1_+)=1$, where $1_+$ is the identity of $A_+$.
Each algebra homomorphism $\varphi\colon A\to B$ uniquely extends
to a unital homomorphism $\varphi_+\colon A_+\to B_+$.
The one-dimensional left (respectively, right) $A$-module $A_+/A$
will be denoted by $\CC_\ell$ (respectively, $\CC_r$).
Note that $\CC_\ell^*$ is isomorphic to $\CC_r$, and vice versa.

The {\em right annihilator} of $A$ is
\[
\rann(A)=\{\, b\in A\; ;\; ab=0\;\forall\, a\in A\}.
\]
If $A$ is a left hypotopological algebra, then it readily follows from the
definition of the first Arens product that
\begin{equation}
\label{rann_bidual}
\rann(A^{**})=\h_A(A^*,\CC_r).
\end{equation}

Let $A$ be a locally convex algebra, and let $A\lcmod$ denote the category
of all left locally convex $A$-modules. Clearly, $A\lcmod$ is an additive category.
For our purposes, it is important that $A\lcmod$ admits inverse and direct limits.
Indeed, if
$\cX=(X_\alpha,\varphi_{\alpha\beta})_{\alpha\in\Lambda}$ is a direct system in
$A\lcmod$ (i.e., a covariant functor from a directed set $\Lambda$ to $A\lcmod$),
then the direct limit of $\cX$ taken in
the category of all left $A$-modules and endowed with the locally convex
direct limit topology is easily seen to be a direct limit of $\cX$
in $A\lcmod$. A similar argument works for inverse limits.

If $A$ is a $\Ptens$-algebra, then the full subcategory of $A\lcmod$
consisting of $\Ptens$-modules will be denoted by $A\lmod$.
Given a full additive subcategory $\cC$ of $\hLCS$, we denote by
$A\lmod(\cC)$ the full subcategory of $A\lmod$ consisting of the modules
whose underlying l.c.s.'s are objects of $\cC$. The respective categories
of right modules and bimodules will be denoted by $\rmod A(\cC)$ and
$A\bimod A(\cC)$ (or just by $\rmod A$ and $A\bimod A$ in the case where
$\cC=\hLCS$). If $A$ is unital, then $A\lunmod(\cC)$ (respectively, $\runmod A(\cC)$,
$A\biunmod A(\cC)$) stands for the full
subcategory of $A\lmod(\cC)$ (respectively, $\rmod A(\cC)$,
$A\bimod A(\cC)$) consisting of unital modules.
Let us remark that $A\lmod(\cC)$ is isomorphic to $A_+\lunmod(\cC)$, and
$A\bimod A(\cC)$ is isomorphic to $A^e\lunmod(\cC)$, where $A^e=A_+\Ptens A_+^\op$
is the enveloping algebra of $A$ \cite[II.5.11]{X1}.
Note also that inverse limits
exist in $A\lmod$ and coincide with inverse limits taken in $A\lcmod$.

If $\cC$ is any of the categories $\hLCS$, $\Fr$, $\Ban$, then each morphism
$\varphi\colon X\to Y$ in $A\lmod(\cC)$ has a kernel and a cokernel.
Specifically, $\Ker\varphi=\varphi^{-1}(0)$ endowed with the
relative topology, and $\Coker\varphi=\bigl(Y/\overline{\varphi(X)}\bigr)\sptilde$,
where $Y/\overline{\varphi(X)}$ is endowed with the quotient topology.
The morphism $\ker\varphi$ is just the embedding of $\varphi^{-1}(0)$
into $X$, and $\coker\varphi$ is the composition of the quotient map
$Y\to Y/\overline{\varphi(X)}$ with the embedding of the target
into its completion. Note that if $\cC$ is either $\Fr$ or $\Ban$, then
$Y/\overline{\varphi(X)}$ is already complete.
A morphism $\varphi\colon X\to Y$ is a kernel
(i.e., there exists a morphism $\psi\colon Y\to Z$ such that
$\varphi=\ker\psi$) if and only if $\varphi$ is {\em topologically injective},
i.e., is a homeomorphism onto $\varphi(X)$.
Similarly, $\varphi\colon X\to Y$
is a cokernel if and only if $\varphi$ is an open map onto a dense submodule of $Y$.
If $\cC$ is either $\Fr$ or $\Ban$, then the Open Mapping Theorem
implies that $\varphi$ is a kernel (respectively, a cokernel) if and only if it is
injective with closed range (respectively, if and only if it is onto).
For $A=0$ (i.e., in the case where $A\lmod(\cC)=\cC$),
the above facts can be found in \cite[Section 3]{Prosm_DFA},
and the same argument works for any $\Ptens$-algebra $A$.

Let $E$ be a vector space, and let $p$ be a seminorm on $E$.
We denote by $E_p$ the completion of $E$ w.r.t. $p$, and
by $\tau_p$ the canonical map $E\to E_p$. Note that $\Ker\tau_p=p^{-1}(0)$.
If $U\subset E$ is an absorbent, absolutely convex set with Minkowski
functional $p_U$, we write $E_U$ for $E_{p_U}$ and $\tau_U$ for $\tau_{p_U}$.
Suppose that $E$ carries a locally convex topology determined by a directed
family of seminorms $\{p_\lambda=\|\cdot\|_\lambda : \lambda\in\Lambda\}$.
For each $\lambda\in\Lambda$, set $E_\lambda=E_{p_\lambda}$
and $\tau_\lambda=\tau_{p_\lambda}$.
Given $\lambda,\mu\in\Lambda$, we write $\lambda\prec\mu$ if
$\| x\|_\lambda\le\| x\|_\mu$
for all $x\in E$.
If $\lambda,\mu\in\Lambda$ and $\lambda\prec\mu$, then
there is a unique continuous linear map $\tau^\mu_\lambda\colon E_\mu\to E_\lambda$
such that $\tau_\lambda=\tau^\mu_\lambda \tau_\mu$. The family
$\{\tau_\lambda : \lambda\in\Lambda\}$ determines a continuous linear map
from $E$ to the inverse limit $\varprojlim (E_\lambda,\tau^\mu_\lambda)$,
which, in turn, yields a topological isomorphism between
$E\sptilde$ and $\varprojlim (E_\lambda,\tau^\mu_\lambda)$
(see, e.g., \cite[II.5.4]{Sch}). Note that the inverse system
$(E_\lambda,\tau^\mu_\lambda)$ is {\em reduced}, i.e., for each
$\mu\in\Lambda$ the canonical map $\varprojlim E_\lambda\to E_\mu$ has dense range.
Conversely, if a complete locally convex space $E$ is isomorphic to a reduced
inverse limit $\varprojlim (E_\lambda,\tau^\mu_\lambda)$ of Banach spaces,
then the topology on $E$ is determined by the family of seminorms
$\{p_\lambda=\|\tau_\lambda(\,\cdot\,)\| : \lambda\in\Lambda\}$,
where $\tau_\lambda\colon E\to E_\lambda$ is the canonical map,
and the inverse system $(E_\lambda,\tau^\mu_\lambda)$ is isomorphic
to the system constructed from the family $\{ p_\lambda\}$ as described above.

If $A$ is an algebra and $p$ is a submultiplicative seminorm on $A$,
then $A_p$ is a Banach algebra and $\tau_p$ is an algebra homomorphism.
Therefore, if $A$ is an Arens-Michael algebra with topology determined
by a directed family $\{\|\cdot\|_\lambda : \lambda\in\Lambda\}$
of submultiplicative seminorms, then all the maps $\tau_\lambda$
and $\tau^\mu_\lambda$ are continuous algebra homomorphisms,
and so we obtain a topological algebra isomorphism
$A\cong\varprojlim A_\lambda$ (the {\em Arens-Michael decomposition theorem}).
The situation described above is
usually expressed by the phrase {\em ``Let $A$ be an Arens-Michael
algebra, and let $A=\varprojlim A_\lambda$ be an Arens-Michael decomposition
of $A$''}.

If $A$ is a locally convex algebra and $X$ is a left locally convex $A$-module,
we say that a continuous seminorm $q$ on $X$ is {\em $m$-compatible}
if there exists a continuous submultiplicative seminorm $p$ on $A$ such that
$q(a\cdot x)\le p(a)q(x)$. In this case, $X_q$ is a Banach $A_p$-module
(and hence a Banach $A$-module) in a natural way, and $\tau_q\colon X\to X_q$
is an $A$-module morphism.
If $A$ is an Arens-Michael algebra and $X$ is a left
$A$-$\Ptens$-module, then the topology on $X$ can be determined by a directed
family $\{\|\cdot\|_\lambda : \lambda\in\Lambda\}$
of $m$-compatible seminorms \cite[3.4]{Pir_qfree}.
Now the maps $\tau_\lambda$
and $\tau^\mu_\lambda$ are morphisms in $A\lmod$,
and so $X$ is isomorphic to $\varprojlim X_\lambda$ in $A\lmod$.
As above, we express this situation by the phrase
{\em ``Let $A$ be an Arens-Michael
algebra, let $X$ be a left $A$-$\Ptens$-module, and let
$X=\varprojlim X_\lambda$ be an Arens-Michael decomposition
of $X$''}.

We end this section by recalling a result of Palamodov on the vanishing
of the derived inverse limit functor
(\cite{Pal_projlim,Pal_methods}; see also \cite[Chapter 3]{Wngnrth}).
The theorem below is a special case of \cite[Corollary 5.1]{Pal_methods};
cf. also \cite[Theorems 4 and 16]{Allan_stabseq}
and \cite[Theorem 2.5]{Allan_auto}.

\begin{theorem}[Palamodov]
\label{thm:Pal}
Let $\cX=(X_n,\varphi^m_n)_{n\in\N}$ be an inverse sequence of Banach spaces
with linking maps $\varphi^m_n\colon X_m\to X_n\; (n\le m)$.
Then the following conditions are equivalent:
\begin{enumerate}
\renewcommand{\theenumi}{\roman{enumi}}
\item
$\varprojlim^1\cX=0$;
\item
for each $n\in\N$ there exists $m\ge n$ such that for each $k\ge m$
the subspace $\varphi^k_n(X_k)$ is dense in $\varphi^m_n(X_m)$.
\end{enumerate}
\end{theorem}

\section{Some properties of the projective tensor product}
\label{sect:tensprod}
Let $A$ be an algebra, $X$ a right $A$-module,
$Y$ a left $A$-module, and $E$ a vector space. Recall that
a bilinear map $\varPhi\colon X\times Y\to E$ is {\em $A$-bilinear}
(or {\em balanced}) if $\varPhi(x\cdot a,y)=\varPhi(x,a\cdot y)$
for all $x\in X,\; y\in Y,\; a\in A$.
If $A$ is a topological algebra, $X,Y$ are topological $A$-modules,
and $E$ is a topological vector space, then the space of all jointly
continuous $A$-bilinear maps from $X\times Y$ to $E$ will be denoted by
$\Bil_A(X\times Y,E)$.

Now suppose that $A$ is a $\Ptens$-algebra, $X,Y$ are $\Ptens$-modules,
and $E$ is a complete l.c.s. Then there exists a natural vector space isomorphism
\[
\Bil_A(X\times Y,E)\cong\cL(X\ptens{A} Y,E),
\]
where $X\ptens{A} Y$ is the {\em $A$-module projective tensor product}
of $X$ and $Y$ \cite[II.4]{X1}.
By definition, $X\ptens{A} Y=(X\Ptens Y/N)\sptilde$, where
$N\subset X\Ptens Y$ is the closure of the linear subspace generated
by all elements of the form $x\cdot a\otimes y-x\otimes a\cdot y\quad
(x\in X,\; y\in Y,\; a\in A)$.

The following two propositions are easy generalizations of the corresponding
results on Banach modules (see \cite[II.3.17]{X1} and
\cite[II.5.21]{X1}, respectively). For Fr\'echet modules,
Proposition~\ref{prop:tens_quot} was also proved in \cite{Podara}.

\begin{prop}
\label{prop:tens_quot}
Let $A$ be a $\Ptens$-algebra, and let $I\subset A_+$ be a closed left ideal.
Then for each $X\in \rmod A$ there is a topological isomorphism
\begin{equation*}
X\ptens{A} (A_+/I)\sptilde\cong (X/\overline{X\cdot I})\sptilde
\end{equation*}
uniquely determined by $x\otimes (a+I)\mapsto x\cdot a+\overline{X\cdot I}$.
\end{prop}
\begin{proof}
For each complete l.c.s. $E$ we have natural isomorphisms
\[
\begin{split}
\cL\bigl(X\ptens{A} (A_+/I)\sptilde,E\bigr)
&\cong \Bil_A(X\times (A_+/I)\sptilde,E)\\
&\cong \Bil_A\bigl(X\times (A_+/I),E\bigr)\\
&\cong \bigl\{\varPhi\in\Bil_A(X\times A_+,E)\, ;\; \varPhi|_{X\times I}=0\bigr\}\\
&\cong \bigl\{\varphi\in\cL(X,E)\, ;\; \varphi|_{X\cdot I}=0\bigr\}\\
&\cong \cL\bigl((X/\overline{X\cdot I})\sptilde,E\bigr).
\end{split}
\]
Setting $E=(X/\overline{X\cdot I})\sptilde$ and applying the above to the identity
map of $E$, we get the result.
\end{proof}

\begin{prop}
\label{prop:adj_ass}
Let $A$ be a $\Ptens$-algebra, $X$ a right $A$-$\Ptens$-module,
and $Y$ a left $A$-$\Ptens$-module. There is a natural linear map
\begin{equation}
\label{adj_ass}
(X\ptens{A} Y)^*\to\h_A(X,Y^*),\qquad
f\mapsto (x\mapsto (y\mapsto f(x\otimes y))).
\end{equation}
The above map is a vector space isomorphism in either of the following cases:
\begin{enumerate}
\renewcommand{\theenumi}{\roman{enumi}}
\item $Y$ is a Banach $A$-module;
\item both $X$ and $Y$ are Fr\'echet $A$-modules.
\end{enumerate}
\end{prop}
\begin{proof}
Clearly, \eqref{adj_ass} is injective. Given
$\varphi\in \h_A(X,Y^*)$, define an $A$-bilinear map $\varPhi\colon X\times Y\to\CC$
by $\varPhi(x,y)=\varphi(x)(y)$.
It is easy to see that $\varPhi$ is separately continuous, and, moreover, is
hypocontinuous w.r.t. the family of all bounded
subsets of $Y$. In either of cases (i) and (ii) this implies that
$\varPhi$ is jointly continuous (case (i) is obvious, while case (ii)
follows from \cite[III.5.1]{Sch}). Therefore there exists $f\in (X\ptens{A} Y)^*$
such that $f(x\otimes y)=\varPhi(x,y)$ for all $x\in X,\; y\in Y$.
Thus \eqref{adj_ass} is onto, as required.
\end{proof}

The projective tensor product preserves cokernels (cf. \cite[V.5.1]{MacLane}).

\begin{prop}
\label{prop:ptenscoker}
Let $A$ be a $\Ptens$-algebra, and let
\[
X_1 \xra{f} X_2 \xra{g} X_3 \to 0
\]
be a sequence in $\rmod A$ such that $g=\coker f$.
Take $Y\in A\lmod$, and consider the sequence
\[
X_1\ptens{A} Y \xra{f\otimes\id_Y} X_2\ptens{A} Y \xra{g\otimes\id_Y}
X_3\ptens{A}Y \to 0.
\]
Then $g\otimes\id_Y=\coker (f\otimes\id_Y)$ in $\hLCS$.
\end{prop}
\begin{proof}
We have to show that for each complete l.c.s. $E$ the sequence
\[
0\to \cL(X_3\ptens{A}Y,E) \to \cL(X_2\ptens{A}Y,E) \to \cL(X_1\ptens{A}Y,E)
\]
is exact in $\Vect$. By the definition of the projective tensor product,
the latter sequence can be identified with
\begin{equation}
\label{bilseq}
0\to \Bil_A(X_3\times Y,E) \xra{\alpha} \Bil_A(X_2\times Y,E) \xra{\beta}
\Bil_A(X_1\times Y,E)
\end{equation}
where $\alpha$ and $\beta$ are induced by $g$ and $f$ in the obvious way.
Since $g(X_2)$ is dense in $X_3$, we see that $\Ker\alpha=0$.
Now suppose that $\varPhi\in\Ker\beta$, i.e., that $\varPhi(x,y)=0$
for all $x\in f(X_1)$ and all $y\in Y$. Since $f(X_1)$ is dense in $\Ker g$,
we have $\varPhi(x,y)=0$ for all $x\in\Ker g$.
Therefore there exists a unique $A$-bilinear map $\varPsi_0\colon g(X_2)\times Y\to E$
such that $\varPsi_0(g(x),y)=\varPhi(x,y)$ for all $x\in X_2,\; y\in Y$.
Since $g\colon X_2\to g(X_2)$ is open, we see that $\varPsi_0$ is jointly
continuous. As $g(X_2)$ is dense in $X_3$, this implies that
$\varPsi_0$ uniquely extends to a map
$\varPsi\in\Bil_A(X_3\times Y,E)$. Clearly, $\alpha(\varPsi)=\varPhi$.
Thus \eqref{bilseq} is exact, and so $g\otimes\id_Y=\coker (f\otimes\id_Y)$,
as required.
\end{proof}

\begin{remark}
In the setting of Banach modules over Banach algebras, Proposition~\ref{prop:ptenscoker}
becomes an easy consequence of \cite[V.5]{ML_work}. Indeed, in this case the
functor $(\,\cdot\, )\ptens{A} Y$ has a right adjoint (namely, $\cL(Y,\,\cdot\,)$)
and therefore preserves colimits.
\end{remark}

\section{Flat and strictly flat Fr\'echet modules}
\label{sect:flat}
Let $\A$ be an additive category with kernels and cokernels.
An {\em exact pair} in $\A$ is a sequence
$X\xra{i} Y\xra{p} Z$ such that $i=\ker p$ and $p=\coker i$.

\begin{remark}
\label{rem:ex_pair}
It is elementary to check that $X\xra{i} Y\xra{p} Z$ is an exact pair
if and only if $p=\coker i$ and $i$ is a kernel, or, equivalently,
if and only if $i=\ker p$ and $p$ is a cokernel.
\end{remark}

By a {\em pre-exact} category we mean an additive category $\A$
endowed with a class $\cE$ of exact pairs closed under isomorphism.
Exact pairs belonging to $\cE$ will be called {\em admissible pairs}
or {\em short admissible sequences}.
A morphism $i\colon X\to Y$ (respectively, a morphism $p\colon Y\to Z$)
is an {\em admissible monomorphism} (respectively, an {\em admissible epimorphism})
if it fits into an admissible pair $X\xra{i} Y\xra{p} Z$.
A chain complex $C_\bullet=(C_n,d_n)_{n\in\Z}$ in $\A$ is {\em admissible}
if for each $n\in\Z$ the morphism $d_n\colon C_n\to C_{n-1}$ has a kernel,
and the sequence $\Ker d_n\to C_n\to\Ker d_{n-1}$ is admissible.

\begin{example}
Each additive category $\A$ becomes a pre-exact category if we define $\cE$ to
be the class of all exact pairs in $\A$. In this case we say that
$\A$ is endowed with the {\em strongest pre-exact structure}.
In the sequel, the categories $\Vect$, $\LCS$, $\hLCS$, $\Fr$, and $\Ban$
will be considered as pre-exact categories with respect to the
strongest pre-exact structure. If $A$ is a $\Ptens$-algebra
and $\cC\subset\hLCS$ is a full additive subcategory, then the category
$A\lmod(\cC)$ endowed with the strongest pre-exact structure will be denoted
by $A\lbarmod(\cC)$. The strongest pre-exact categories $\rbarmod A(\cC)$
and $A\bibarmod A(\cC)$
(as well as $A\lunbarmod(\cC)$, $\runbarmod A(\cC)$ and $A\biunbarmod A(\cC)$
for $A$ unital) are defined similarly. Note that if $\cC$ is either $\Ban$ or $\Fr$,
then a chain complex in any of the above categories is admissible if and
only if it is (algebraically) exact.
\end{example}

\begin{remark}
As was observed in \cite{Prosm_DFA}, the categories $\cC=\LCS$, $\cC=\Fr$,
and $\cC=\Ban$
are in fact exact categories in the sense of Quillen \cite{Quillen},
but this is not the case for $\cC=\hLCS$. Similar statements hold for
the respective $A$-module categories $A\lbarmod(\cC)$,
$\rbarmod A(\cC)$, $A\bibarmod A(\cC)$.
This is the reason why we have to consider pre-exact categories that are not exact.
\end{remark}

\begin{example}
There is another way to make a category of topological modules
into a pre-exact category.
Let $A$ be a $\Ptens$-algebra, and let $\cC$ be a full additive subcategory
of $\hLCS$. We say that an exact pair $X\to Y\to Z$
in $A\lmod(\cC)$ is admissible if the sequence $0\to X\to Y\to Z\to 0$
splits in $\cC$. The resulting pre-exact category will be denoted by
the same symbol $A\lmod(\cC)$. The pre-exact categories $\rmod A(\cC)$
and $A\bimod A(\cC)$ are defined similarly. These are the basic
categories of Topological Homology (see \cite{X1}). In fact, it is easy to
show that the above categories are exact in the sense of Quillen.
\end{example}

Let $\A$ and $\B$ be pre-exact categories, and let
$F\colon \A\to\B$ be an additive functor. We say that $F$ is {\em exact}
if for each admissible pair $X\to Y\to Z$ in $\A$ the pair
$FX\to FY\to FZ$ is admissible in $\B$. A standard argument shows that
$F$ is exact if and only if for each admissible complex $C_\bullet$
in $\A$ the complex $F(C_\bullet)$ is admissible in $\B$.

\begin{definition}[Helemskii \cite{Hel_period,Hel_flat_amen,X1}]
\label{def:flat}
Let $A$ be a Fr\'echet algebra.
A left Fr\'echet $A$-module $X$ is {\em flat}
if the functor
\begin{equation}
\label{tens_func}
(\,\cdot\,)\ptens{A} X\colon \rmod A(\Fr)\to\Vect
\end{equation}
is exact. The module $X$ is {\em strictly flat}
if the functor
\begin{equation}
\label{tens_func2}
(\,\cdot\,)\ptens{A} X\colon \rbarmod A(\Fr)\to\Vect
\end{equation}
is exact.
A right Fr\'echet $A$-module (respectively, a Fr\'echet $A$-bimodule)
is flat if it is flat as a left Fr\'echet module over $A^\op$
(respectively, over $A^e$). Strictly flat right Fr\'echet modules
and bimodules are defined similarly.
\end{definition}

\begin{remark}
\label{rem:flat_Fr}
As easily follows from the Open Mapping Theorem, a pair $X\to Y\to Z$
in $\Fr$ is exact if and only if it is exact in $\Vect$.
Therefore we may consider the functors appearing in Definition~\ref{def:flat}
as functors with values in $\Fr$. This will lead to an equivalent definition.
\end{remark}

\begin{remark}
\label{rem:uni_flat}
Suppose that $A$ is a unital Fr\'echet algebra and $X$ is a unital left Fr\'echet $A$-module.
Then $X$ is flat (respectively, strictly flat) if and only if
the restriction of \eqref{tens_func} (respectively, of \eqref{tens_func2})
to $\runmod A(\Fr)$ (respectively, to $\runbarmod A(\Fr)$) is exact.
Indeed, each right Fr\'echet $A$-module $Y$ can be decomposed as
$Y=Y^\un\oplus Y^\ann$, where $Y^\un$ is a unital $A$-module and $Y^\ann$
is an annihilator $A$-module (i.e., $y\cdot a=0$ for all $y\in Y^\ann,\; a\in A$).
Specifically, $Y^\un=\{ y\in Y\, :\; y\cdot 1=y\}$ and
$Y^\ann=\{ y\in Y\, :\; y\cdot 1=0\}$. It is easy to see that this decomposition
is natural, so that each complex $Y_\bullet$ of right Fr\'echet $A$-modules decomposes
as $Y_\bullet=Y_\bullet^\un\oplus Y_\bullet^\ann$.
Clearly, $Y_\bullet$ is admissible in $\rmod A(\Fr)$ (respectively, in
$\rbarmod A(\Fr)$) if and only if both $Y^\un_\bullet$ and $Y^\ann_\bullet$ are.
On the other hand, if $X$ is unital, then $Y_\bullet^\ann\ptens{A} X=0$,
and so $Y_\bullet\ptens{A} X\cong Y^\un_\bullet\ptens{A} X$.
This proves the claim.
\end{remark}

\begin{remark}
Certainly, Definition~\ref{def:flat} makes sense for more general module
categories of the form $A\lmod(\cC)$ and $A\lbarmod(\cC)$. However,
for the reasons explained in Remark~\ref{rem:flat_LCS2} below, we shall consider
flatness and strict flatness only in categories of Fr\'echet modules
over Fr\'echet algebras.
\end{remark}

\begin{prop}
\label{prop:flat_base}
Let $\varphi\colon A\to B$ be a Fr\'echet algebra homomorphism, and let
$X$ be a flat (respectively, strictly flat) left Fr\'echet
$A$-module. Then $B_+\ptens{A} X$ is flat (respectively, strictly flat)
in $B\lmod(\Fr)$.
\end{prop}
\begin{proof}
The functors $(\,\cdot\,)\ptens{B} B_+\ptens{A} X$ and
$(\,\cdot\,)\ptens{A} X$ on $\rmod B(\Fr)$ are obviously isomorphic.
\end{proof}

\begin{prop}
\label{prop:flat-3}
Let $A$ be a Fr\'echet algebra, and let
\[
0 \to X\to Y\to Z\to 0
\]
be an admissible sequence in $A\lmod(\Fr)$. Suppose that $Z$ is flat.
Then $X$ is flat if and only if $Y$ is flat.
\end{prop}
\begin{proof}
This is a special case of \cite[Proposition 3.1.15]{EschmPut}.
\end{proof}

Flat modules are traditionally considered within the category of Banach
modules over a Banach algebra.
The definition of a flat Banach module, as given in \cite{Hel_period,Hel_flat_amen,X1},
is formally different from Definition~\ref{def:flat}, as the domain of the
functor $(\,\cdot\,)\ptens{A} X$ is taken to be $\rmod A(\Ban)$
rather than $\rmod A(\Fr)$. The following proposition shows that the two definitions
of flatness are consistent.

\begin{prop}
\label{prop:flat_Ban_fr}
Let $A$ be a Banach algebra, and let $X$ be a left Banach $A$-module.
Then the following conditions are equivalent:
\begin{enumerate}
\renewcommand{\theenumi}{\roman{enumi}}
\item the functor
\[
(\,\cdot\,)\ptens{A} X\colon \rmod A(\Fr)\to\Vect
\]
is exact (i.e., $X$ is flat in the sense of Definition~\ref{def:flat});
\item the functor
\[
(\,\cdot\,)\ptens{A} X\colon \rmod A(\Ban)\to\Vect
\]
is exact (i.e., $X$ is flat in the sense of \cite{Hel_period}).
\end{enumerate}
\end{prop}
\begin{proof}
Clearly, (i) implies (ii). Conversely, suppose that (ii) holds, and let
$Y_\bullet$ be an admissible complex in $\rmod A(\Fr)$. In order to prove that
$Y_\bullet\ptens{A} X$ is exact, it suffices to prove that
$(Y_\bullet\ptens{A} X)^*$ is exact (see, e.g., \cite[26.4]{MV}).
By Proposition~\ref{prop:adj_ass}, we have an isomorphism
\[
(Y_\bullet\ptens{A} X)^* \cong \h_A(Y_\bullet,X^*).
\]
Since $X^*$ is injective in $\rmod A(\Ban)$ \cite[VII.1.14]{X1}, it follows that
$X^*$ is a retract of $\cL(A_+,X^*)$ in $\rmod A(\Ban)$
\cite[III.1.31]{X1}. On the other hand,
\[
\h_A(Y_\bullet,\cL(A_+,X^*))\cong \cL(Y_\bullet,X^*)
\]
(see \cite[Proposition 3.2]{T1}), and the latter complex is obviously exact.
This completes the proof.
\end{proof}

A similar statement holds for strictly flat Banach modules.

\begin{prop}
Let $A$ be a Banach algebra, and let $X$ be a left Banach $A$-module.
Then the following conditions are equivalent:
\begin{enumerate}
\renewcommand{\theenumi}{\roman{enumi}}
\item the functor
\[
(\,\cdot\,)\ptens{A} X\colon \rbarmod A(\Fr)\to\Vect
\]
is exact (i.e., $X$ is strictly flat in the sense of Definition~\ref{def:flat});
\item the functor
\[
(\,\cdot\,)\ptens{A} X\colon \rbarmod A(\Ban)\to\Vect
\]
is exact (i.e., $X$ is strictly flat in the sense of \cite{Hel_flat_amen}).
\end{enumerate}
\end{prop}
\begin{proof}
Suppose that (ii) holds, and let $Y_\bullet$ be an exact complex in
$\rbarmod A(\Fr)$. As in Proposition~\ref{prop:flat_Ban_fr}, we have
to prove that $\h_A(Y_\bullet,X^*)$ is exact.
By \cite[VII.1.14]{X1}, $X^*$ is strictly injective in $\rmod A(\Ban)$.
Choose an injective Banach space $E$ of the form $E=\ell^\infty(S)$
and an isometric embedding
$\varphi\colon X^*\to E$. Define an $A$-module morphism $\varkappa$ to be the composition
\[
X^* \xra{\nu} \cL(A_+,X^*) \xra{\varphi_*} \cL(A_+,E),
\]
where $\nu$ is the canonical embedding of $X^*$ into $\cL(A_+,X^*)$
given by $\nu(f)(a)=f\cdot a\; (f\in X^*,\; a\in A_+)$
(see \cite[Section III.1.4]{X1}), and $\varphi_*$
is induced by $\varphi$ in the obvious way. Since both $\nu$ and $\varphi_*$
are topologically injective, so is $\varkappa=\varphi_*\nu$.
Since $X^*$ is strictly injective in $\rmod A(\Ban)$, there exists
an $A$-module morphism $\sigma\colon\cL(A_+,E)\to X^*$ such that
$\sigma\varkappa=\id_{X^*}$ \cite[VII.1.13]{X1}.
Therefore $X^*$ is a retract of $\cL(A_+,E)$
in $\rmod A(\Ban)$. On the other hand,
\[
\h_A(Y_\bullet,\cL(A_+,E))\cong \cL(Y_\bullet,E)
\]
(see \cite[Proposition 3.2]{T1}), and the latter complex is exact
due to the injectivity of $E=\ell^\infty(S)$
in $\Fr$ (see, e.g., \cite[\S4]{Pal_methods} or \cite[2.2.1]{Wngnrth}).
This completes the proof.
\end{proof}

Recall that a Banach algebra $A$ is {\em amenable} \cite{Jhnsn_CBA} if for
each Banach $A$-bimodule $X$ every continuous derivation from $A$
to the dual bimodule, $X^*$, is inner. As was shown by Helemskii and Sheinberg
\cite{Hel_Shein} (see also \cite[VII.2.17]{X1}), $A$ is amenable
if and only if $A_+$ is a flat Banach $A$-bimodule.

By definition \cite{X_31,X_HOA},
a Fr\'echet algebra $A$ is {\em amenable} (respectively, {\em biflat})
if $A_+$ (respectively, $A$) is a flat Fr\'echet $A$-bimodule.
A number of interesting characterizations of biflat Banach algebras was
given by Selivanov \cite{Sel_cohchar}.

Using Proposition~\ref{prop:flat_Ban_fr}, we get the following.

\begin{corollary}
If $A$ is an amenable (respectively, biflat) Banach algebra, then $A$
is amenable (respectively, biflat) when considered as a Fr\'echet algebra.
\end{corollary}

\begin{prop}
\label{prop:amen_dense}
Let $\varphi\colon A\to B$ be a Fr\'echet algebra homomorphism.
Suppose that $A$ is amenable and that the map
\begin{equation}
\label{bprod}
B_+\ptens{A_+} B_+\to B_+,\qquad b_1\otimes b_2\mapsto b_1 b_2
\end{equation}
is bijective. Then $B$ is amenable.
\end{prop}
\begin{proof}
We have
\[
B^e\ptens{A^e} A_+\cong B_+\ptens{A_+} A_+\ptens{A_+} B_+
\cong B_+\ptens{A_+} B_+\cong B_+,
\]
so the result follows from Proposition~\ref{prop:flat_base}.
\end{proof}

\begin{remark}
\label{rem:amen_dense}
It can be shown that \eqref{bprod} is bijective if and only if
$\varphi$ is an epimorphism in the category of Fr\'echet algebras
(cf. \cite[XI.1.2]{Stnstrm}). Note that this condition is satisfied
whenever $\varphi$ has dense range.
\end{remark}

Given a Fr\'echet algebra $A$ and a Fr\'echet $A$-bimodule $X$, consider
the chain complex
\begin{equation}
\label{std_chain}
0 \lar C_0(A,X) \xla{d_1} C_1(A,X) \xla{d_2} \cdots
\xla{d_n} C_n(A,X) \xla{d_{n+1}} \cdots
\end{equation}
where $C_0(A,X)=X, \; C_n(A,X)=X\Ptens\underbrace{A\Ptens\cdots\Ptens A}_n$
for $n\ge 1$, and the differentials are given by
\[
\begin{split}
d_n(x\otimes a_1\otimes\cdots\otimes a_n)
&=x\cdot a_1\otimes a_2\otimes\cdots\otimes a_n\\
&+\sum_{k=1}^{n-1} (-1)^k x\otimes a_1\otimes\cdots\otimes a_k a_{k+1}
\otimes\cdots\otimes a_n\\
&+(-1)^n a_n\cdot x\otimes a_1\otimes\cdots\otimes a_{n-1},\\
d_1(x\otimes a)&=x\cdot a-a\cdot x.
\end{split}
\]
By definition \cite{Guichardet1} (see also \cite{X1,Dales}),
the $n$th homology of \eqref{std_chain}
is denoted by $\cH_n(A,X)$ and is called the {\em $n$th Hochschild homology group
of A with coefficients in $X$}.

By \cite[VII.2.17]{X1}, $A$ is amenable if and only if for each $X\in A\bimod A(\Fr)$
the space $\cH_0(A,X)$ is Hausdorff, and $\cH_n(A,X)=0$ for all $n\ge 1$.
Denote by $\overline{[A,X]}\subset X$ the closure of the image of $d_1$.
Then it follows from the above that $A$ is amenable if and only if
the augmented complex
\begin{equation}
\label{std_hom_augm}
0 \lar X/\overline{[A,X]} \lar C_\bullet(A,X)
\end{equation}
is exact.

\section{Approximate identities and approximate diagonals\\
in topological algebras}
\label{sect:bai}

In this section we collect some general facts on approximate identities
and approximate diagonals. Although
these facts are well known in the Banach algebra case (see, e.g., \cite{Dales,Palmer}),
we give full proofs in most cases for the sake of completeness.

Recall that a net $(e_\alpha)_{\alpha\in\Lambda}$ in a topological algebra $A$
is a {\em right approximate identity} (a right a.i. for short)
if $a=\lim_\alpha ae_\alpha$ for each $a\in A$.
Left approximate identities are defined similarly.
A net $(e_\alpha)_{\alpha\in\Lambda}$ is a {\em two-sided approximate identity}
(or just an {\em approximate identity}) if it is both a left and a right
approximate identity. An approximate identity
$(e_\alpha)_{\alpha\in\Lambda}$ (right, left, or two-sided)
is {\em bounded} if the set
$\{ e_\alpha\}_{\alpha\in\Lambda}$ is bounded.
The expression ``bounded approximate identity'' is abbreviated as ``b.a.i.'',
as usual.

We start with a simple observation.
Let us recall that a Hausdorff locally convex space $E$ is {\em semi-Montel} if each bounded
subset of $E$ is relatively compact.

\begin{lemma}
\label{lemma:bai_Montel}
Let $A$ be a semi-Montel topological algebra with a right b.a.i.
Then $A$ has a right identity.
\end{lemma}
\begin{proof}
Let $(e_\alpha)_{\alpha\in\Lambda}$ be a right b.a.i. in $A$.
By assumption, the set $\overline{ \{ e_\alpha\}_{\alpha\in\Lambda} }$
is compact, so that there exists a subnet $(e_{\alpha'})_{\alpha'\in\Lambda'}$
of $(e_\alpha)$ convergent to $e\in A$. Then for each $a\in A$ the net
$(ae_{\alpha'})$ converges to $a$ and to $ae$ simultaneously.
Hence $a=ae$, and so $e$ is a right identity in $A$.
\end{proof}

\begin{prop}
\label{prop:ai_crit}
Let $A$ be a topological algebra. Then
\begin{enumerate}
\renewcommand{\theenumi}{\roman{enumi}}
\item $A$ has a right a.i. if and only if for each finite subset $F\subset A$
and each $0$-neighborhood $U\subset A$ there exists $b\in A$ such that
$a-ab\in U$ for all $a\in F$;
\item $A$ has a right b.a.i. if and only if there exists a bounded subset $B\subset A$
such that for each finite subset $F\subset A$
and each $0$-neighborhood $U\subset A$ there exists $b\in B$ such that
$a-ab\in U$ for all $a\in F$.
\end{enumerate}
\end{prop}
\begin{proof}
If $(e_\alpha)$ is a right a.i. (respectively, a right b.a.i.) in $A$,
then $a-ae_\alpha\to 0$ uniformly
on finite subsets of $A$. Therefore the element $b=e_\alpha$
satisfies (i) (respectively, (ii)) for $\alpha$ large enough.
To prove the converse, consider the set $S$ of all pairs $(U,F)$, where $U\subset A$
is a $0$-neighborhood and $F\subset A$ is a finite set, and define an order
on $S$ by
\[
(U_1,F_1)\le (U_2,F_2) \quad\iff\quad U_1\supset U_2\quad\text{and}\quad F_1\subset F_2.
\]
Thus $S$ becomes a directed set. For each $\alpha=(U,F)\in S$ find $e_\alpha\in A$
(respectively, $e_\alpha\in B$)
such that $a-ae_\alpha\in U$ for all $a\in F$. Then it is clear that $(e_\alpha)$
is a right a.i.
(respectively, a right b.a.i.) in $A$.
\end{proof}

\begin{remark}
\label{rem:ai_crit}
Proposition~\ref{prop:ai_crit} has obvious ``left'' and ``two-sided'' versions.
In particular, $A$ has a two-sided (bounded) a.i.
if and only if (there exists a bounded subset $B\subset A$ such that)
for each finite subset $F\subset A$
and each $0$-neighborhood $U\subset A$ there exists $b\in A$ (respectively, $b\in B$)
such that $a-ab\in U$ and $a-ba\in U$ for all $a\in F$.
\end{remark}

\begin{remark}
\label{rem:lrbai2bai}
As in the case of Banach algebras, it can be shown that if a
hypotopological algebra has a left and a right b.a.i., then it has a two-sided b.a.i.
We omit the proof, as we do not use this result in the sequel.
\end{remark}

\begin{prop}
\label{prop:bai_id_bidual}
Let $A$ be a left hypotopological algebra. The following conditions are equivalent:
\begin{enumerate}
\renewcommand{\theenumi}{\roman{enumi}}
\item
$A$ has a right b.a.i.;
\item
there exists a bounded set $B\subset A$ such that $a\in \overline{aB}$
for each $a\in A$.
\end{enumerate}
If, in addition, $A$ is locally convex, then {\upshape (i)} and {\upshape (ii)}
are equivalent to
\begin{enumerate}
\item[{\upshape (iii)}]
$A^{**}$ has a right identity;
\item[{\upshape (iv)}]
there exists $e\in A^{**}$ such that $i_A(a)e=i_A(a)$ for each $a\in A$.
\end{enumerate}
\end{prop}
\begin{proof}
$\mathrm{(i)}\Longrightarrow\mathrm{(ii)}$. This is clear.

$\mathrm{(ii)}\Longrightarrow\mathrm{(i)}$.
Let $F=\{ a_1,\ldots ,a_n\}\subset A$ be a finite set and $U\subset A$ a $0$-neighborhood.
Choose a circled $0$-neighborhood $V\subset A$ such that $V+V+V\subset U$
and a $0$-neighborhood $W\subset V$ such that $WB\subset V$ and $FW\subset V$.
For each $j=1,\ldots ,n$ set $B_j=(1_+-B)^{n-j}\subset A_+$.
Since $A$ (and hence $A_+$) is left hypotopological, we see that $B_j$ is bounded
(cf. Section~\ref{sect:prelim} or \cite{Jhnsn_centr}).
Let $W^+\subset A_+$ be a $0$-neighborhood such that $W^+\cap A=W$.
For each $j=1,\ldots ,n$ choose a $0$-neighborhood
$W_j^+\subset A_+$ such that $W_j^+ B_j\subset W^+$. Setting $W_j=W_j^+\cap A$,
we see that $W_j B_j\subset W$.

Using (ii), we can inductively find $u_1,\ldots ,u_n\in B$ such that
\[
a_j(1_+-u_1)\cdots (1_+-u_j)\in W_j\qquad (j=1,\ldots ,n).
\]
Define $v\in A$ by $1_+-v=(1_+-u_1)\cdots (1_+-u_n)$. Then
\[
\begin{split}
a_j(1_+-v)=a_j(1_+-u_1)\cdots (1_+-u_j)&(1_+-u_{j+1})\cdots (1_+-u_n)\\
&\in W_j B_j\subset W
\qquad (j=1,\ldots ,n).
\end{split}
\]
Finally, choose $u\in B$ such that $v-vu\in W$. Then for each $j=1,\ldots ,n$ we have
\[
\begin{split}
a_j-a_ju&=(a_j-a_jv)+a_j(v-vu)-(a_j-a_jv)u\\
&\in W+FW-WB\subset V+V+V\subset U.
\end{split}
\]
Applying Proposition~\ref{prop:ai_crit}, we see that $A$ has a right b.a.i.

Now suppose that $A$ is locally convex.

$\mathrm{(i)}\Longrightarrow\mathrm{(iv)}$.
Let $(e_\alpha)_{\alpha\in\Lambda}$ be a right b.a.i. in $A$.
Since $\{ i_A(e_\alpha)\}\subset A^{**}$ is equicontinuous, there exists
a subnet $(e_{\alpha'})_{\alpha'\in\Lambda'}$ of $(e_\alpha)_{\alpha\in\Lambda}$
such that $(i_A(e_{\alpha'}))$
weak$^*$ converges to $e\in A^{**}$. This implies that
$(i_A(a)i_A(e_{\alpha'}))$ weak$^*$ converges to $i_A(a)e$ for each $a\in A$
(\cite[3.6]{Gulick}; see also Section~\ref{sect:prelim}).
On the other hand, $ae_{\alpha'}\to a$, so that $i_A(a)i_A(e_{\alpha'})\to i_A(a)$,
and, finally, $i_A(a)=i_A(a)e$ for each $a\in A$.

$\mathrm{(iv)}\Longrightarrow\mathrm{(iii)}$.
This follows form the weak$^*$ density of $i_A(A)$ in $A^{**}$ and
the weak$^*$ continuity of the first Arens product in the first variable
(\cite[3.4]{Gulick}).

$\mathrm{(iii)}\Longrightarrow\mathrm{(ii)}$.
Let $e\in A^{**}$ be a right identity. By the bipolar theorem, there exists
a convex bounded subset $B\subset A$ such that $e$ belongs to the weak$^*$ closure
of $i_A(B)$. Using again \cite[3.6]{Gulick}, we see that for each $a\in A$
the element $i_A(a)e=i_A(a)$ belongs to the weak$^*$ closure of $i_A(a)i_A(B)=i_A(aB)$.
Equivalently, $a$ belongs to the weak closure of $aB$.
Since $aB$ is convex, this completes the proof of (ii).
\end{proof}

Let $A$ be a $\Ptens$-algebra, and let $\pi_0\colon A\Ptens A\to A$
denote the product map. (We write $\pi_0$ instead of $\pi$ in order
to make our notation consistent with that of \cite[VII.2.2]{X1}
and of Section~\ref{sect:flat_amen}).
Following \cite{Jhnsn_appr}, we say that a net
$(M_\alpha)_{\alpha\in\Lambda}$ in $A\Ptens A$
is an {\em approximate diagonal} for $A$
if $a\cdot M_\alpha-M_\alpha\cdot a\to 0$ and
$\pi_0(M_\alpha)a\to a$ for each $a\in A$.
An approximate diagonal $(M_\alpha)_{\alpha\in\Lambda}$ is bounded
if the set $\{ M_\alpha\}_{\alpha\in\Lambda}$ is bounded.
An element $M\in (A\Ptens A)^{**}$ is a {\em virtual diagonal} for $A$
if $a\cdot M=M\cdot a$ and $\pi_0^{**}(M)\cdot a=i_A(a)$
for each $a\in A$. (Here $(A\Ptens A)^{**}$ is considered as
a semitopological $A$-bimodule; see Section~\ref{sect:prelim}.)

In what follows, given $a\in A$ and $M\in A\Ptens A$, we set
$[a,M]=a\cdot M-M\cdot a$.

\begin{remark}
Our terminology slightly differs from that of \cite{Jhnsn_appr}.
Namely, our ``bounded approximate diagonal'' is just ``approximate diagonal''
in \cite{Jhnsn_appr}. However, it has recently become clear that
unbounded approximate diagonals are also of interest
(cf. \cite{GL,DLZ}).
\end{remark}

\begin{lemma}
\label{lemma:bad_vd}
Let $A$ be a $\Ptens$-algebra. Then $A$ has a bounded approximate diagonal
if and only if $A$ has a virtual diagonal.
\end{lemma}

We omit the proof, as it is similar to the Banach algebra case \cite[1.2]{Jhnsn_appr}
(see also \cite[2.9.64]{Dales}).

Let us recall some standard notation (see, e.g., \cite{Sch}).
Suppose that $\mathscr U$ is a neighborhood base at $0$ in $A$.
Given $U\in\mathscr U$, let $\overline{\Gamma(U\otimes U)}$ denote the
closure of the absolutely convex hull of the set
\[
U\otimes U=\{ a\otimes b : a,b\in U\}\subset A\Ptens A.
\]
Then $\{ \overline{\Gamma(U\otimes U)} : U\in\mathscr U\}$
is a neighborhood base at $0$ in $A\Ptens A$ \cite{Sch}.

\begin{prop}
\label{prop:ad_crit}
Let $A$ be a $\Ptens$-algebra. Then
\begin{enumerate}
\renewcommand{\theenumi}{\roman{enumi}}
\item $A$ has an approximate diagonal if and only if for each
finite subset $F\subset A$
and each $0$-neighborhood $U\subset A$ there exists $M\in A\Ptens A$ with
$[a,M]\in\overline{\Gamma(U\otimes U)}$
and $\pi_0(M)a-a\in U$ for all $a\in F$;
\item $A$ has a bounded approximate diagonal if and only if
there exists a bounded subset $B\subset A\Ptens A$ such that
for each finite subset $F\subset A$
and each $0$-neighborhood $U\subset A$ there exists $M\in B$ with
$[a,M]\in\overline{\Gamma(U\otimes U)}$
and $\pi_0(M)a-a\in U$ for all $a\in F$.
\end{enumerate}
\end{prop}

The proof is similar to that of Proposition~\ref{prop:ai_crit} and is therefore
omitted.

\section{Locally bounded approximate identities and\\ locally bounded
approximate diagonals}
\label{sect:lbai}
\begin{definition}
\label{def:lbai}
Let $A$ be a topological algebra. We say that
\begin{enumerate}
\renewcommand{\theenumi}{\roman{enumi}}
\item $A$ {\em has a right (respectively, left)
locally bounded approximate identity}
if for each $0$-neighborhood
$U\subset A$ there exists $C>0$ such that for each
finite subset $F\subset A$ there exists $b\in CU$
with $a-ab\in U$ (respectively, $a-ba\in U$) for all $a\in F$;
\item $A$ {\em has a two-sided locally bounded approximate identity}
(or just a {\em locally bounded approximate identity})
if for each $0$-neighborhood
$U\subset A$ there exists $C>0$ such that for each
finite subset $F\subset A$ there exists $b\in CU$
with $a-ab\in U$ and $a-ba\in U$ for all $a\in F$.
\end{enumerate}
We abbreviate ``locally bounded approximate identity'' as ``locally b.a.i.''
\end{definition}

\begin{remark}
\label{rem:def_base}
A slightly more convenient definition is as follows: $A$ has a right
locally b.a.i. if for each $0$-neighborhood
$U\subset A$ there exists $C>0$ such that for each
finite subset $F\subset A$ and each $\eps>0$ there exists $b\in CU$
with $a-ab\in \eps U$ for all $a\in F$.
To show that these definitions are equivalent, it suffices to
substitute $\eps^{-1}F$ for $F$ in Definition~\ref{def:lbai}.
It is also clear that instead of considering arbitrary $0$-neighborhoods
in $A$, we can restrict to $0$-neighborhoods belonging to some family $\mathscr U$
such that $\{ tU : t>0,\; U\in\mathscr U\}$ is a neighborhood base at $0$.
\end{remark}

\begin{remark}
\label{rem:lbai_semi}
Let $A$ be a locally convex algebra, and let
$\{\|\cdot\|_\lambda : \lambda\in\Lambda\}$ be a directed family
of seminorms generating the topology of $A$. Using Remark~\ref{rem:def_base},
we see that $A$ has a right locally b.a.i.
if and only if there exists a family $\{ C_\lambda : \lambda\in\Lambda\}$
of positive reals such that for each finite subset $F\subset A$,
each $\lambda\in\Lambda$, and each $\eps>0$ there exists $b\in A$
with $\| b\|_\lambda\le C_\lambda$ and
$\| a-ab\|_\lambda<\eps$
for all $a\in F$.
Similar statements hold for left and for
two-sided locally b.a.i.'s.

Let us remark that the above characterization
of locally b.a.i.'s in terms of seminorms was taken as definition
in \cite{Pir_ICTAA05}. We thank A.~Mallios for suggesting us to give a definition
applicable to arbitrary (not necessarily locally convex) topological algebras
(see Definition~\ref{def:lbai}).
\end{remark}

\begin{remark}
Taking into account Proposition~\ref{prop:ai_crit}, we see that the existence
of a right (respectively, left, two-sided) b.a.i. implies the existence of
a right (respectively, left, two-sided) locally b.a.i., which, in turn, implies
the existence of a right (respectively, left, two-sided) a.i.
\end{remark}

\begin{remark}
\label{rem:normable}
If $A$ is normable, then the notions of ``bounded'' and ``locally bounded''
a.i.'s are obviously equivalent (cf. Remark~\ref{rem:lbai_semi}).
\end{remark}

\begin{prop}
\label{prop:dense}
Let $\varphi\colon A\to B$ be a continuous homomorphism
of topological algebras. Suppose that $A$ has a right locally b.a.i.
Then each of the following conditions implies that $B$ has a right
locally b.a.i.:
\begin{enumerate}
\renewcommand{\theenumi}{\roman{enumi}}
\item $\varphi$ is onto;
\item $\varphi$ has dense range, and the multiplication in $B$ is jointly
continuous.
\end{enumerate}
\end{prop}
\begin{proof}
We prove (ii), the proof of (i) being similar.
Given a $0$-neighborhood $V\subset B$, find circled $0$-neighborhoods
$V_1\subset V$ and $V_2\subset V_1$ such that $V_1+V_1+V_1\subset V$
and $V_2^2\subset V_1$. Set $U=\varphi^{-1}(V_2)$,
and let $C>0$ be as in Definition~\ref{def:lbai}.
Without loss of generality, we assume that $C\ge 1$.
Given a finite set $F'\subset B$, find a finite set $F\subset A$
such that $F'\subset\varphi(F)+C^{-1} V_2$. Then there exists $b\in CU$
such that $a-ab\in U$ for each $a\in F$. Take any $a'\in F'$ and choose $a\in F$
satisfying $a'-\varphi(a)\in C^{-1} V_2$.
Then $b'=\varphi(b)\in CV_2\subset CV$, and
\[
\begin{split}
a'-a'b'&=\varphi(a-ab)+(a'-\varphi(a))+(\varphi(a)-a')b'\\
&\in V_1+ V_1+ C^{-1}V_2\cdot CV_2
\subset V_1+ V_1+ V_1\subset V.
\end{split}
\]
The rest is clear.
\end{proof}

\begin{prop}
\label{prop:dense2}
Let $A$ be an Arens-Michael algebra. Then the following conditions are equivalent:
\begin{enumerate}
\renewcommand{\theenumi}{\roman{enumi}}
\item $A$ has a right locally b.a.i.;
\item for each Banach algebra $B$ such that there exists a continuous homomorphism
$\varphi\colon A\to B$ with dense range, the algebra $B$ has a right b.a.i.
\end{enumerate}
\end{prop}
\begin{proof}
$\mathrm{(i)}\Longrightarrow\mathrm{(ii)}$.
This follows from Proposition~\ref{prop:dense} and Remark~\ref{rem:normable}.

$\mathrm{(ii)}\Longrightarrow\mathrm{(i)}$.
Let $U\subset A$ be an absolutely convex, idempotent $0$-neighborhood.
By assumption, the Banach algebra $A_U$ has a right b.a.i.
Hence the same is true of the dense subalgebra $A_U^0=\Im\tau_U\subset A_U$
(see \cite[1.4]{DW}). Therefore there exists $C>0$ such that for each finite set
$F\subset A$ there exists $b\in A$ such that
$\|\tau_U(b)\|<C$ and that $\| \tau_U(a)-\tau_U(a)\tau_U(b)\|<1$
for each $a\in F$. This means exactly that $b\in CU$ and that
$a-ab\in U$ for each $a\in F$. In view of Remark~\ref{rem:def_base},
this completes the proof.
\end{proof}

\begin{corollary}
\label{cor:bai_AMdec}
Let $A$ be a Arens-Michael algebra, and let $A=\varprojlim A_\lambda$ be an
Arens-Michael decomposition of $A$. Then $A$ has a right locally b.a.i. if
and only if each $A_\lambda$ has a right b.a.i.
\end{corollary}
\begin{proof}
Each continuous homomorphism $\varphi$ from $A$ to a Banach
algebra $B$ uniquely factors through some $A_\lambda$.
Therefore if $\varphi$ has dense range and $A_\lambda$ has a right
b.a.i., then the same is true of $B$ (cf. \cite[1.4]{DW}).
\end{proof}

\begin{remark}
\label{rem:bai_AMdec}
It is easy to see that Propositions~\ref{prop:dense} and ~\ref{prop:dense2}
and Corollary~\ref{cor:bai_AMdec} have obvious ``right'' and ``two-sided''
versions.
\end{remark}

\begin{remark}
As in the case of Banach algebras, it can be shown that if a topological
algebra with jointly continuous multiplication has a left and a right
locally b.a.i., then it has a two-sided locally b.a.i.
(cf. also Remark~\ref{rem:lrbai2bai}).
We omit the proof, as we do not use this result in the sequel.
\end{remark}

\begin{definition}
\label{def:lbad}
Let $A$ be a $\Ptens$-algebra. We say that
$A$ {\em has a locally bounded approximate diagonal}
if for each $0$-neighborhood
$U\subset A$ there exists $C>0$ such that for each
finite subset $F\subset A$ there exists $M\in C\overline{\Gamma(U\otimes U)}$
with $[a,M]\in \overline{\Gamma(U\otimes U)}$ and
$\pi_0(M)a-a\in U$ for all $a\in F$.
\end{definition}

\begin{remark}
Taking into account Proposition~\ref{prop:ad_crit}, we see that the existence
of a bounded approximate diagonal implies the existence of
a locally bounded approximate diagonal, which, in turn, implies
the existence of an approximate diagonal.
\end{remark}

\begin{remark}
\label{rem:normable_diag}
For a Banach algebra, the notions of ``bounded'' and ``locally bounded''
approximate diagonals are equivalent by Proposition~\ref{prop:ad_crit} (ii).
\end{remark}

\begin{prop}
\label{prop:dense_diag}
Let $\varphi\colon A\to B$ be a continuous homomorphism
of $\Ptens$-algebras. Suppose that $A$ has a locally bounded
approximate diagonal and that $\varphi$ has dense range.
Then $B$ has a locally bounded approximate diagonal.
\end{prop}
\begin{proof}
Given an absolutely convex $0$-neighborhood $V\subset B$, find an absolutely convex
$0$-neighborhood
$V_1\subset B$ such that $V_1\cap V_1^2\cap \overline{V_1^3}\subset V/3$.
Set $U=\varphi^{-1}(V_1)$, and let $C>0$ be as in Definition~\ref{def:lbad}.
Without loss of generality, we assume that $C\ge 1$.
Given a finite set $F'\subset B$, find a finite set $F\subset A$
such that $F'\subset\varphi(F)+C^{-1} V_1$. Then there exists
$M\in C\overline{\Gamma(U\otimes U)}$
with $[a,M]\in \overline{\Gamma(U\otimes U)}$ and
$\pi_0(M)a-a\in U$ for all $a\in F$.

Set $M'=(\varphi\otimes\varphi)(M)$. Clearly, $M'\in C\overline{\Gamma(V_1\otimes V_1)}
\subset C\overline{\Gamma(V\otimes V)}$.
Now take any $a'\in F'$ and choose $a\in F$
satisfying $a'-\varphi(a)\in C^{-1} V_1$.
We have
\[
\begin{split}
[a',M']&=(\varphi\otimes\varphi)[a,M]+[a'-\varphi(a),M']\\
&\in \overline{\Gamma(V_1\otimes V_1)}
+ C^{-1} V_1\cdot C\overline{\Gamma(V_1\otimes V_1)}
+C\overline{\Gamma(V_1\otimes V_1)}\cdot C^{-1}V_1\\
&\subset \overline{\Gamma(V\otimes V)}.
\end{split}
\]
Furthermore,
\[
\begin{split}
\pi_0(M')a'-a'&=\varphi(\pi_0(M)a-a)+(\varphi(a)-a')+\pi_0(M')(a'-\varphi(a))\\
&\in V_1+ V_1+ C\overline{V_1^2} C^{-1}V_1\subset V.
\end{split}
\]
The rest is clear.
\end{proof}

\begin{prop}
\label{prop:dense_diag2}
Let $A$ be an Arens-Michael algebra. Then the following conditions are equivalent:
\begin{enumerate}
\renewcommand{\theenumi}{\roman{enumi}}
\item $A$ has a locally bounded approximate diagonal;
\item for each Banach algebra $B$ such that there exists a continuous homomorphism
$\varphi\colon A\to B$ with dense range, the algebra $B$ has a bounded approximate diagonal.
\end{enumerate}
\end{prop}
\begin{proof}
$\mathrm{(i)}\Longrightarrow\mathrm{(ii)}$.
This follows from Proposition~\ref{prop:dense_diag} and Remark~\ref{rem:normable_diag}.

$\mathrm{(ii)}\Longrightarrow\mathrm{(i)}$.
Let $U\subset A$ be an absolutely convex, idempotent $0$-neighborhood.
By assumption, the Banach algebra $A_U$ has a bounded approximate diagonal.
Therefore there exists $C>0$ such that for each finite set
$F\subset A$ there exists $M'\in A_U\Ptens A_U$ with
\begin{align}
\label{diag_const}
&\| M'\|<C,\\
\label{diag_commut}
&\|\, [\tau_U(a),M']\,\| < 1\quad \forall a\in F,\\
\label{diag_ai}
&\| \pi_0(M')\tau_U(a)-\tau_U(a)\| < 1\quad \forall a\in F.
\end{align}
Clearly, the set of all $M'\in A_U\Ptens A_U$ satisfying
\eqref{diag_const}--\eqref{diag_ai} is open. On the other hand,
$(\tau_U\otimes\tau_U)(A)$ is dense in $A_U\Ptens A_U$,
and so there exists $M\in A\Ptens A$ such that $M'=(\tau_U\otimes\tau_U)(M)$
satisfies \eqref{diag_const}--\eqref{diag_ai}.
This implies that $M$ satisfies the conditions of Definition~\ref{def:lbad}.
\end{proof}

\begin{corollary}
\label{cor:bad_AMdec}
Let $A$ be a Arens-Michael algebra, and let $A=\varprojlim A_\lambda$ be an
Arens-Michael decomposition of $A$. Then $A$ has a
locally bounded approximate diagonal if
and only if each $A_\lambda$ has a bounded approximate diagonal.
\end{corollary}

\section{Strictly flat cyclic modules}

In this section we extend the first part of Helemskii and Sheinberg's
Theorem~\ref{thm:Hel_flat} (specifically, the equivalence of conditions (i) and (ii))
to Arens-Michael algebras.

\label{sect:strflat}
\begin{lemma}
\label{lemma:projquot}
Let $A$ be an Arens-Michael algebra, let $X\in \rmod A$, and let
$X=\varprojlim(X_\lambda,\tau^\mu_\lambda)$ be an Arens-Michael decomposition
of $X$. Suppose that $I\subset A$ is a closed left ideal.
For each $\lambda$ consider the map
\[
\hat\tau_\lambda\colon X/\overline{X\cdot I}\to X_\lambda/\overline{X_\lambda\cdot I},
\qquad
x+\overline{X\cdot I}\mapsto \tau_\lambda(x)+\overline{X_\lambda\cdot I}.
\]
Then $(X_\lambda/\overline{X_\lambda\cdot I},\hat\tau_\lambda)$ is the completion
of $X/\overline{X\cdot I}$ with respect to the quotient seminorm of $\|\cdot\|_\lambda$.
As a consequence, we have a topological isomorphism
\[
(X/\overline{X\cdot I})\sptilde\cong\varprojlim
(X_\lambda/\overline{X_\lambda\cdot I},\hat\tau^\mu_\lambda),
\]
where the structure maps $\hat\tau^\mu_\lambda\; (\lambda\prec\mu)$ are defined by
\[
\hat\tau^\mu_\lambda\colon X_\mu/\overline{X_\mu\cdot I}
\to X_\lambda/\overline{X_\lambda\cdot I},
\qquad
x+\overline{X_\mu\cdot I}\mapsto \tau^\mu_\lambda(x)+\overline{X_\lambda\cdot I}.
\]
\end{lemma}
\begin{proof}
Since $\tau_\lambda\colon X\to X_\lambda$ has dense range, the same is true of
$\hat\tau_\lambda$. Therefore we need only prove that $\hat\tau_\lambda$ is
isometric with respect to the quotient seminorm of $\|\cdot\|_\lambda$
on $X/\overline{X\cdot I}$ and the quotient norm on
$X_\lambda/\overline{X_\lambda\cdot I}$. Denote the former seminorm by
$\|\cdot\|'_\lambda$ and the latter norm by $\|\cdot\|'$.
Then for each $x\in X$ we have
\[
\begin{split}
\| x+\overline{X\cdot I}\|'_\lambda
&=\inf\{ \| x+y\|_\lambda : y\in \overline{X\cdot I}\}
=\inf\{ \| \tau_\lambda(x)+\tau_\lambda(y)\| : y\in \overline{X\cdot I}\}\\
&=\inf\{ \| \tau_\lambda(x)+z\| : z\in \overline{X_\lambda\cdot I}\}
=\| \tau_\lambda(x)+\overline{X_\lambda\cdot I}\|'\\
&=\| \hat\tau_\lambda(x+\overline{X\cdot I})\|',
\end{split}
\]
which proves the claim.
\end{proof}

\begin{theorem}
\label{thm:strflat}
Let $A$ be an Arens-Michael algebra and $I\subset A_+$ a closed left ideal.
Set $X=(A_+/I)\sptilde$. Then the following conditions are equivalent:
\begin{enumerate}
\renewcommand{\theenumi}{\roman{enumi}}
\item $I$ has a right locally b.a.i.;
\item the functor
$(\,\cdot\,)\ptens{A} X\colon \rbarmod A\to\hLCS$
is exact;
\item the functor
$(\,\cdot\,)\ptens{A} X\colon \rbarmod A(\Ban)\to\Ban$
is exact.
\end{enumerate}
\end{theorem}
\begin{proof}
$\mathrm{(i)}\Longrightarrow\mathrm{(ii)}$.
In view of Proposition~\ref{prop:ptenscoker} and Remark~\ref{rem:ex_pair},
it suffices to prove that
for each $Y\in\rbarmod A$ and each closed submodule $Z\subset Y$
the map
\begin{equation}
\label{tens_id}
Z\ptens{A} X\to Y\ptens{A} X
\end{equation}
is topologically injective. Choose an Arens-Michael decomposition
$Y=\varprojlim Y_\lambda$ of $Y$, and, for each $\lambda$, set
$Z_\lambda=\overline{\tau_\lambda(Z)}$, where $\tau_\lambda\colon Y\to Y_\lambda$
is the canonical map.
Then it is easy to see that $Z=\varprojlim Z_\lambda$ is an Arens-Michael
decomposition of $Z$ (cf. also \cite[2.5.6]{Eng}).
Since each seminorm $\|\cdot\|_\lambda$ on $Y$ is
$m$-compatible (see Section~\ref{sect:prelim}), there exists a continuous
submultiplicative seminorm $\|\cdot\|_\lambda^\circ$ on $A$ such that
$\| y\cdot a\|_\lambda\le \| y\|_\lambda \| a\|_\lambda^\circ$ for each
$y\in Y,\; a\in A$. Let $A_\lambda$ be the completion of $A$ with respect to
$\|\cdot\|_\lambda^\circ$. Then it is clear that $Y_\lambda$ is a right Banach
$A_\lambda$-module in a canonical way, and that $Z_\lambda$ is a closed
$A_\lambda$-submodule of $X_\lambda$.

Let $\sigma_\lambda\colon A\to A_\lambda$ be the canonical map, and let
$I_\lambda=\overline{(\sigma_\lambda)_+(I)}$. Clearly, $I_\lambda$ is a closed
left ideal in $(A_\lambda)_+$. By Proposition~\ref{prop:dense2}, $I_\lambda$
has a right b.a.i., and so $X_\lambda=(A_\lambda)_+/I_\lambda$ is a strictly
flat Banach $A_\lambda$-module (see Theorem~\ref{thm:Hel_flat}). Therefore the map
\[
Z_\lambda\ptens{A_\lambda} X_\lambda\to Y_\lambda\ptens{A_\lambda} X_\lambda
\]
is topologically injective. Using the identifications
\[
\xymatrix@R-5pt{
Z_\lambda\ptens{A_\lambda} X_\lambda \ar[r] \ar@{=}[d]_{\textstyle\wr}
& Y_\lambda\ptens{A_\lambda} X_\lambda \ar@{=}[d]^{\textstyle\wr} \\
Z_\lambda/\overline{Z_\lambda\cdot I_\lambda} \ar[r] \ar@{=}[d]
& Y_\lambda/\overline{Y_\lambda\cdot I_\lambda} \ar@{=}[d]\\
Z_\lambda/\overline{Z_\lambda\cdot I} \ar[r]
& Y_\lambda/\overline{Y_\lambda\cdot I}
}
\]
and taking the inverse limit, we get a topologically injective map
\[
\varprojlim Z_\lambda/\overline{Z_\lambda\cdot I}
\to\varprojlim Y_\lambda/\overline{Y_\lambda\cdot I}.
\]
On the other hand, Lemma~\ref{lemma:projquot} and Proposition~\ref{prop:tens_quot}
imply that the latter map is identifiable with \eqref{tens_id} as follows:
\[
\xymatrix@R-5pt{
\varprojlim Z_\lambda/\overline{Z_\lambda\cdot I} \ar[r] \ar@{=}[d]_{\textstyle\wr}
& \varprojlim Y_\lambda/\overline{Y_\lambda\cdot I} \ar@{=}[d]^{\textstyle\wr} \\
(Z/\overline{Z\cdot I})\sptilde \ar[r] \ar@{=}[d]_{\textstyle\wr}
& (Y/\overline{Y\cdot I})\sptilde \ar@{=}[d]^{\textstyle\wr}\\
Z\ptens{A} X \ar[r]
& Y\ptens{A} X.
}
\]
This proves (ii).

$\mathrm{(ii)}\Longrightarrow\mathrm{(iii)}$. This is clear.

$\mathrm{(iii)}\Longrightarrow\mathrm{(i)}$.
Let $\varphi\colon A\to B$ be a homomorphism of $A$ to a Banach algebra $B$.
Using the same argument as in Proposition~\ref{prop:flat_base},
we see that the left Banach $B$-module $B_+\ptens{A} X\cong B_+/\overline{B_+\cdot I}$
is strictly flat. Therefore the ideal $\overline{B_+\cdot I}$ of $B_+$
has a right b.a.i.

Now choose an Arens-Michael decomposition $A=\varprojlim A_\lambda$
and apply the above result to the canonical homomorphism
$\sigma_\lambda\colon A\to A_\lambda$. In this case, the ideal
$I_\lambda=\overline{(A_\lambda)_+\cdot I}$ of $(A_\lambda)_+$
is equal to the closure of $(\sigma_\lambda)_+(I)$ in $(A_\lambda)_+$,
so that $I=\varprojlim I_\lambda$ is an Arens-Michael decomposition of $I$.
By Corollary~\ref{cor:bai_AMdec}, $I$ has a right locally b.a.i.,
as required.
\end{proof}

\begin{corollary}
\label{cor:strflat}
Let $A$ be a Fr\'echet-Arens-Michael algebra
and $I\subset A_+$ a closed left ideal.
Set $X=A_+/I$. Then the following conditions are equivalent:
\begin{enumerate}
\renewcommand{\theenumi}{\roman{enumi}}
\item $I$ has a right locally b.a.i.;
\item $X$ is a strictly flat Fr\'echet $A$-module.
\end{enumerate}
\end{corollary}

\begin{remark}
\label{rem:strflat_unital}
Theorem~\ref{thm:strflat} and Corollary~\ref{cor:strflat}
have obvious ``unital'' versions. Specifically, if $A$ is unital, then
we can replace $A_+$ by $A$ and ``$\underline{\mathsf{mod}}$''
by ``$\underline{\mathsf{unmod}}$'' everywhere in the statements of
Theorem~\ref{thm:strflat} and Corollary~\ref{cor:strflat}
(cf. also Remark~\ref{rem:uni_flat}).
\end{remark}

\begin{remark}
\label{rem:flat_LCS2}
In view of Theorem~\ref{thm:strflat} and Corollary~\ref{cor:strflat},
it is tempting to call a left $\Ptens$-module $X$ over a $\Ptens$-algebra
$A$ strictly flat if it satisfies condition (ii) of Theorem~\ref{thm:strflat}.
However, it is not clear whether such a ``strictly flat'' module is always flat
(which, by definition \cite{X1}, means that the functor
$(\,\cdot\,)\ptens{A} X\colon\rmod A\to\Vect$ is exact). A possible way out of
this situation is to redefine the notion of flatness as follows:
{\em $X\in A\lmod(\cC)$ is flat if the functor
$(\,\cdot\,)\ptens{A} X\colon \rmod A(\cC)\to\cC$ is exact}.
If $\cC=\Fr$ or $\cC=\Ban$, this is equivalent to Definition~\ref{def:flat}
(cf. Remark~\ref{rem:flat_Fr}).
However, the disadvantage of the new definition is that it might be
impossible to characterize such ``flatness'' in terms of the derived functors
$\Tor_n$ (cf. \cite{X1}).
That is why we do not consider flat and strictly flat $\Ptens$-modules in the
general setting, restricting ourselves to the metrizable case.
\end{remark}

\section{Quasinormable Fr\'echet algebras}
\label{sect:qsinorm}
Recall that a Hausdorff l.c.s. $E$ is {\em quasinormable} \cite{Groth_F_DF}
if for each $0$-neighborhood $U\subset E$ there exists a $0$-neighborhood
$V\subset U$ such that for each $\eps>0$ there exists a bounded set $B\subset E$
satisfying $V\subset B+\eps U$. In this section we prove that
for quasinormable Fr\'echet algebras
the notions of bounded and locally bounded approximate identities
are equivalent.

\begin{lemma}
\label{lemma:ran_bidual}
Let $\varphi\colon A\to B$ be a Banach algebra homomorphism with dense range.
Then
\begin{enumerate}
\renewcommand{\theenumi}{\roman{enumi}}
\item $\varphi^{**}(\rann(A^{**}))\subset \rann(B^{**})$;
\item if $e$ is a right identity in $A^{**}$, then
$\varphi^{**}(e)$ is a right identity in $B^{**}$.
\end{enumerate}
\end{lemma}
\begin{proof}
If $a\in \rann(A^{**})$, then $b\varphi^{**}(a)=0$ for all $b\in\Im\varphi^{**}$.
Since $\Im\varphi$ is dense in $B$, it follows that $\Im\varphi^{**}$ is
weak$^*$ dense in $B^{**}$. As the product in $B^{**}$
is weak$^*$ continuous w.r.t. the first variable, this implies that
$b\varphi^{**}(a)=0$ for all $b\in B^{**}$, i.e., that
$\varphi^{**}(a)\in \rann(B^{**})$. Similarly, if $e$ is a right identity
in $A^{**}$, then $\varphi^{**}(e)$ is a right identity for $\Im\varphi^{**}$,
and, by weak$^*$ continuity, for $B^{**}$.
\end{proof}

From now on, let $A$ be a Fr\'echet-Arens-Michael algebra, and let
$A=\varprojlim(A_n,\tau^m_n)_{n\in\N}$ be an Arens-Michael decomposition of $A$.
By Lemma~\ref{lemma:ran_bidual}, we have a well-defined inverse system
$(\rann(A_n^{**}),\sigma^m_n)$, where $\sigma^m_n$ is the restriction of
$(\tau^m_n)^{**}$ to $\rann(A_m^{**})$.

\begin{lemma}
\label{lemma:ran_lim1}
Suppose that $A$ is quasinormable and has a right locally b.a.i.
Then $\varprojlim^1 \rann(A_n^{**})=0$.
\end{lemma}
\begin{proof}
Since $A$ is quasinormable, the direct system $(A_n^*,(\tau^m_n)^*)$
is acyclic \cite[Prop.~7.5]{Pal_methods}, and so $\varprojlim^1 A_n^{**}=0$
\cite[Prop.~6.2]{Pal_methods}. By Palamodov's Theorem~\ref{thm:Pal},
for each $n\in\N$ there exists $m\ge n$ such that
$(\tau^k_n)^{**}(A_k^{**})$
is dense in $(\tau^m_n)^{**}(A_m^{**})$ for all $k\ge m$.

For each $n\in\N$ set $R_n=\rann(A_n^{**})$.
We claim that $(R_n,\sigma^m_n)$ satisfies the conditions
of Palamodov's Theorem with the same $m=m(n)$. To prove the claim,
fix $n\in\N$, and take any $k\ge m$. Since $A$ has a right locally b.a.i.,
Corollary~\ref{cor:bai_AMdec} implies that $A_k$ has a right b.a.i.,
and so $A_k^{**}$ has a right identity, say $e_k$. By Lemma~\ref{lemma:ran_bidual},
$e_m=(\tau^k_m)^{**}(e_k)$ and $e_n=(\tau^m_n)^{**}(e_m)$ are right identities
in $A_m^{**}$ and $A_n^{**}$, respectively.

Define $p_k\colon A_k^{**}\to R_k$ by $p_k(a)=a-e_k a$.
Clearly, $p_k$ is a projection onto $R_k$.
Similarly, define projections
$p_m\colon A_m^{**}\to R_m,\; a\mapsto a-e_m a$,
and $p_n\colon A_n^{**}\to R_n,\; a\mapsto a-e_n a$.
We obtain a commutative diagram
\[
\xymatrix@C+15pt{
A_n^{**} \ar[d]_{p_n} & A_m^{**} \ar[l]^{(\tau^m_n)^{**}} \ar[d]_{p_m}
& A_k^{**} \ar[l]^{(\tau^k_m)^{**}} \ar[d]^{p_k} \\
R_n & R_m \ar[l]^{\sigma^m_n}
& R_k \ar[l]^{\sigma^k_m}
}
\]
Therefore,
\[
\begin{split}
\overline{\sigma^m_n(R_m)}=\overline{\sigma^m_n(p_m(A_m^{**}))}
&=\overline{p_n(\overline{(\tau^m_n)^{**}(A_m^{**})})}\\
&=\overline{p_n(\overline{(\tau^k_n)^{**}(A_k^{**})})}
=\overline{\sigma^k_n(p_k(A_k^{**}))}
=\overline{\sigma^k_n(R_k)}.
\end{split}
\]
Thus $(\rann(A_n^{**}),\sigma^m_n)$ satisfies the conditions
of Palamodov's Theorem \ref{thm:Pal}, and so
$\varprojlim^1 \rann(A_n^{**})=0$.
\end{proof}

\begin{remark}
We note that Lemma~\ref{lemma:ran_lim1} holds under a weaker assumption
than the quasinormability of $A$. Specifically, it suffices to require that
the direct system $(A_n^*,(\tau^m_n)^*)$ be {\em weakly acyclic}
\cite{Pal_methods}; see also \cite[Chapter 6]{Wngnrth}.
\end{remark}

\begin{theorem}
\label{thm:quasinorm}
Let $A$ be a quasinormable Fr\'echet-Arens-Michael algebra with a right
locally b.a.i. Then $A$ has a right b.a.i.
\end{theorem}
\begin{proof}
Let $A=\varprojlim(A_n,\tau^m_n)_{n\in\N}$ be an Arens-Michael
decomposition of $A$. For each $n\in\N$ consider the admissible sequences
\begin{gather*}
0 \lar \CC_\ell \xla{\eps_{A_n}} (A_n)_+ \lar A_n \lar 0,\\
0 \to \CC_r \xra{\eps_{A_n}^*} (A_n)_+^* \to A_n^* \to 0
\end{gather*}
in $A_n\lmod(\Ban)$ and in $\rmod A_n(\Ban)$, respectively.
Since $A$ has a right locally b.a.i.,
Corollary~\ref{cor:bai_AMdec} implies that $A_n$ has a right b.a.i.,
and so $\CC_r$ is injective in $\rmod A_n(\Ban)$ \cite[VII.1.20]{X1}.
Hence the sequence
\begin{equation}
\label{hom1}
0 \to \h_{A_n} (A_n^*,\CC_r) \to \h_{A_n} ((A_n)_+^*,\CC_r)
\to \h_{A_n} (\CC_r,\CC_r) \to 0
\end{equation}
is exact. Since $\tau_n\colon A\to A_n$ has dense range,
we have $\h_{A_n}=\h_A$ in the above sequence. Using the isomorphism
$\h_{A_n} (\CC_r,\CC_r)\cong\CC$, we may identify \eqref{hom1}
with the sequence
\begin{equation}
\label{hom2}
0 \to \h_A (A_n^*,\CC_r) \to \h_A ((A_n)_+^*,\CC_r)
\xra{\varPhi_n} \CC \to 0,
\end{equation}
where $\varPhi_n(f)=f(\eps_{A_n})$ for each $f\in \h_A ((A_n)_+^*,\CC_r)$.

Using \eqref{rann_bidual}, we see that
$\h_A (A_n^*,\CC_r)=\h_{A_n} (A_n^*,\CC_r)=\rann(A_n^{**})$,
and it is clear that for each $m\ge n$ the map
\[
\h_A (A_m^*,\CC_r)\to \h_A (A_n^*,\CC_r),\quad f\mapsto f\circ (\tau^m_n)^*
\]
coincides with the map $\sigma^m_n$ defined before Lemma~\ref{lemma:ran_lim1}.
Therefore taking the inverse limit in \eqref{hom2} and applying
Lemma~\ref{lemma:ran_lim1} we get an exact sequence
\begin{equation}
\label{hom3}
0 \to \varprojlim\h_A (A_n^*,\CC_r) \to \varprojlim\h_A ((A_n)_+^*,\CC_r)
\to \CC \to 0.
\end{equation}
Since $A$ is quasinormable (hence distinguished \cite[26.19]{MV}),
the canonical linear map $\varinjlim A_n^*\to A^*$ is an
isomorphism of locally convex spaces \cite[25.13]{MV}.
Clearly, this is also an isomorphism in the category of right locally convex $A$-modules.
Hence we obtain vector space isomorphisms
\[
\varprojlim\h_A (A_n^*,\CC_r)
\cong \h_A(\varinjlim A_n^*,\CC_r)
\cong \h_A(A^*,\CC_r).
\]
Similarly, $\varprojlim\h_A ((A_n)_+^*,\CC_r)\cong \h_A(A_+^*,\CC_r)$ as vector spaces.
Therefore \eqref{hom3} can be identified with
\begin{equation}
\label{hom4}
0 \to \h_A (A^*,\CC_r) \to \h_A (A_+^*,\CC_r)
\xra{\varPhi} \CC \to 0,
\end{equation}
where $\varPhi(f)=f(\eps_A)$ for each $f\in \h_A (A_+^*,\CC_r)$.

As \eqref{hom4} is exact, there exists $f\in \h_A (A_+^*,\CC_r)$
such that $f(\eps_A)=1$. This is equivalent to say that
the sequence
\[
0 \to \CC_r \xra{\eps_A^*} A_+^* \xra{j_A^*} A^* \to 0
\]
splits in the category of right locally convex $A$-modules.
Therefore the dual sequence,
\[
0 \lar \CC_\ell \xla{\eps_A^{**}} A_+^{**} \xla{j_A^{**}} A^{**} \lar 0\, ,
\]
splits in $A\lmod(\Fr)$, and so there exists a left Fr\'echet $A$-module
morphism $\sigma\colon A_+^{**}\to A^{**}$ such that $\sigma j_A^{**}=\id_{A^{**}}$.
Set $e=\sigma(i_{A_+}(1_+))\in A^{**}$. For each $a\in A$ we have
\[
i_A(a)e=a\cdot e=\sigma(i_{A_+}(a\cdot 1_+))=\sigma(i_{A_+}(j_A(a))
=\sigma(j_A^{**}(i_A(a))=i_A(a).
\]
By Proposition~\ref{prop:bai_id_bidual}, this is equivalent to say that
$A$ has a right b.a.i., as required.
\end{proof}

\section{Flat Fr\'echet modules and amenable Fr\'echet algebras}
\label{sect:flat_amen}
Part (iii) of Helemskii and Sheinberg's Theorem~\ref{thm:Hel_flat} can be generalized
to Arens-Michael algebras as follows.

\begin{theorem}
\label{thm:flat}
Let $A$ be an Arens-Michael algebra and $I\subset A_+$ a closed left ideal.
Set $X=(A_+/I)\sptilde$. Suppose that there exists an Arens-Michael decomposition
$A=\varprojlim (A_\lambda,\sigma^\lambda_\mu)_{\lambda\in\Lambda}$
such that for each $\lambda\in\Lambda$ the ideal
$I_\lambda=\overline{(\sigma_\lambda)_+(I)}$ is weakly complemented in $(A_\lambda)_+$.
Then the following conditions are equivalent:
\begin{enumerate}
\renewcommand{\theenumi}{\roman{enumi}}
\item $I$ has a right locally b.a.i.;
\item the functor
$(\,\cdot\,)\ptens{A} X\colon \rbarmod A\to\hLCS$
is exact;
\item the functor
$(\,\cdot\,)\ptens{A} X\colon \rbarmod A(\Ban)\to\Ban$
is exact;
\item the functor
$(\,\cdot\,)\ptens{A} X\colon \rmod A\to\hLCS$
is exact;
\item the functor
$(\,\cdot\,)\ptens{A} X\colon \rmod A(\Ban)\to\Ban$
is exact.
\end{enumerate}
\end{theorem}
\begin{proof}
$\mathrm{(i)}\iff\mathrm{(ii)}\iff\mathrm{(iii)}$: this is Theorem~\ref{thm:strflat}.

$\mathrm{(ii)}\Longrightarrow \mathrm{(iv)}\Longrightarrow \mathrm{(v)}$: this is clear.

$\mathrm{(v)}\Longrightarrow \mathrm{(i)}$. The proof is similar to the proof of
implication $\mathrm{(iii)}\Longrightarrow \mathrm{(i)}$ in Theorem~\ref{thm:strflat};
one uses the equivalence of (i) and (iii) in Theorem~\ref{thm:Hel_flat}.
\end{proof}

\begin{corollary}
\label{cor:flat}
Let $A$ be a Fr\'echet-Arens-Michael algebra
and $I\subset A_+$ a closed left ideal.
Set $X=A_+/I$.
Suppose that there exists an Arens-Michael decomposition
$A=\varprojlim (A_\lambda,\sigma^\lambda_\mu)_{\lambda\in\N}$
such that for each $\lambda\in\N$ the ideal
$I_\lambda=\overline{(\sigma_\lambda)_+(I)}$ is weakly complemented in $(A_\lambda)_+$.
Then the following conditions are equivalent:
\begin{enumerate}
\renewcommand{\theenumi}{\roman{enumi}}
\item $I$ has a right locally b.a.i.;
\item $X$ is a strictly flat Fr\'echet $A$-module;
\item $X$ is a flat Fr\'echet $A$-module.
\end{enumerate}
\end{corollary}

\begin{remark}
\label{rem:flat_unital}
Like Theorem~\ref{thm:strflat} and Corollary~\ref{cor:strflat},
Theorem~\ref{thm:flat} and Corollary~\ref{cor:flat}
have obvious ``unital'' versions. Specifically, if $A$ is unital, then
we can replace $A_+$ by $A$, ``$\underline{\mathsf{mod}}$''
by ``$\underline{\mathsf{unmod}}$'', and
``$\mathsf{mod}$'' by ``$\mathsf{unmod}$'' everywhere in the statements of
Theorem~\ref{thm:flat} and Corollary~\ref{cor:flat}.
\end{remark}

\begin{remark}
The additional condition imposed on $I_\lambda$ in
Theorem~\ref{thm:flat} and in Corollary~\ref{cor:flat}
looks rather unnatural,
but we do not know how to put it into a more reasonable form. In particular,
it is not clear whether or not this condition depends on the particular
choice of an Arens-Michael decomposition. Also we do not know which relation
(if any) this condition has to the position of $I$ inside $A_+$.
Fortunately, as we shall see in Theorem~\ref{thm:amen} below, there is an
important situation where this condition is satisfied automatically.
\end{remark}

Let us recall some notation from \cite[VII.2]{X1}.
Let $A$ be an Arens-Michael algebra, and let $A^e=A_+\Ptens A_+^\op$ be the
enveloping algebra of $A$. Denote by $\pi\colon A^e\to A_+$ the product map,
and set $I^\Delta=\Ker\pi$. Note that $I^\Delta$ is a closed left ideal of
$A^e$, and that $A_+$ is isomorphic to $A^e/I^\Delta$ in $A^e\lunmod=A\bimod A$.
Moreover, $I^\Delta$ is complemented in $A^e$, since
the map $A_+\to A^e,\; a\mapsto 1_+\otimes a$ is readily seen to be a right
inverse of $\pi$.

Similarly, set $A_0^e=A\Ptens A^\op$,
denote by $\pi_0\colon A^e_0\to A$ the product map,
and set $I^\Delta_0=\Ker\pi_0$. Finally, set
$M=\Ker(\eps_A\otimes \eps_A)\subset A^e$.
Note that $M=A_+\Ptens A+A\Ptens A_+$.

Now choose an Arens-Michael decomposition
$A=\varprojlim (A_\lambda,\sigma^\lambda_\mu)_{\lambda\in\Lambda}$.
For each $\lambda\in\Lambda$ we use the symbols $A_\lambda^e,\; (A_0^e)_\lambda,\;
\pi_\lambda,\; (\pi_0)_\lambda,\; I^\Delta_\lambda,\; (I^\Delta_0)_\lambda,\;
M_\lambda$ to denote the respective algebras, maps and ideals constructed
from $A_\lambda$, as described above. For each $\lambda\prec\mu$ we set
\begin{gather*}
\tau^\mu_\lambda=(\sigma^\mu_\lambda)_+\otimes(\sigma^\mu_\lambda)_+\colon
A_\mu^e\to A_\lambda^e,\\
\tau_\lambda=(\sigma_\lambda)_+\otimes(\sigma_\lambda)_+\colon
A^e\to A_\lambda^e.
\end{gather*}

\begin{lemma}
\label{lemma:diagonal}
For each $\lambda\in\Lambda$ we have
\begin{equation}
\label{dense_ideals}
I^\Delta_\lambda=\overline{\tau_\lambda(I^\Delta)},\quad
M_\lambda=\overline{\tau_\lambda(M)},
\end{equation}
and there exist Arens-Michael decompositions
\begin{gather}
\label{AM_dec_tens}
A^e\cong\varprojlim A^e_\lambda, \qquad
A^e_0\cong\varprojlim (A^e_0)_\lambda,\\
\label{AM_dec_diag}
I^\Delta\cong\varprojlim I^\Delta_\lambda,\qquad
M\cong\varprojlim M_\lambda,
\end{gather}
where the linking maps are the restrictions of $\tau^\mu_\lambda$
to the respective subalgebras of $A_\mu^e$.
If, in addition, $A$ has a one-sided locally b.a.i., then
\begin{equation}
\label{dense_ideals2}
(I^\Delta_0)_\lambda=\overline{(\sigma_\lambda\otimes \sigma_\lambda) (I^\Delta_0)},
\end{equation}
and there exists an Arens-Michael decomposition
\begin{equation}
\label{AM_dec_diag2}
I^\Delta_0\cong\varprojlim (I^\Delta_0)_\lambda.
\end{equation}
\end{lemma}
\begin{proof}
Since $M_\lambda=(A_\lambda)_+\Ptens A_\lambda+A_\lambda\Ptens (A_\lambda)_+$,
the second formula in \eqref{dense_ideals} is immediate from
the fact that $\sigma_\lambda\colon A\to A_\lambda$ has dense range.
The first formula in \eqref{dense_ideals} is proved similarly, by using
the fact that $I^\Delta_\lambda$ is the smallest closed
left ideal of $A_\lambda^e$ containing all elements of the form
$1_+\otimes a-a\otimes 1_+\; (a\in (A_\lambda)_+)$ (see \cite[VII.2.11]{X1}).

Now suppose that $A$ has a one-sided locally b.a.i.
By Corollary~\ref{cor:bai_AMdec}, each $A_\lambda$ has a one-sided
b.a.i., and so the product map $A_\lambda\ptens{A_\lambda} A_\lambda\to A_\lambda$
is an isomorphism (see \cite[II.3.12]{X1}). Therefore $(I^\Delta_0)_\lambda$ coincides with
the kernel of the canonical map
\[
\hat\pi_\lambda\colon
A_\lambda\Ptens A_\lambda\to A_\lambda\ptens{A_\lambda} A_\lambda,\quad
a\otimes b\mapsto a\otimes b.
\]
On the other hand, $\Ker\hat\pi_\lambda$ is equal to the closure of the linear span
of all elements of the form $ab\otimes c-a\otimes bc\; (a,b,c\in A_\lambda)$;
see \cite[II.3]{X1}. Using the fact that
$\sigma_\lambda\colon A\to A_\lambda$ has dense range,
we obtain \eqref{dense_ideals2}.

The isomorphisms \eqref{AM_dec_tens} follow from the fact that
projective tensor product commutes with reduced inverse limits \cite[41.6]{Kothe2}.
Now \eqref{AM_dec_diag} and \eqref{AM_dec_diag2}
follow from \eqref{AM_dec_tens}, \eqref{dense_ideals}, and \eqref{dense_ideals2}
(cf. also \cite[2.5.6]{Eng}).
\end{proof}

\begin{lemma}
\label{lemma:bifl_amen}
Let $A$ be a Fr\'echet-Arens-Michael-algebra with a locally b.a.i.
Then $A$ is amenable if and only if $A$ is biflat.
\end{lemma}
\begin{proof}
Let $A=\varprojlim A_\lambda$ be an Arens-Michael decomposition of $A$.
Since $A$ has a locally b.a.i., Corollary~\ref{cor:bai_AMdec} shows that
each $A_\lambda$ has a b.a.i.
By \cite[Lemma VII.2.12 (V)]{X1},
this implies that for each $\lambda\in\N$ the ideal
$M_\lambda$ has a b.a.i. Applying Lemma~\ref{lemma:diagonal}
and Corollary~\ref{cor:bai_AMdec}, we see that $M$ has a locally b.a.i.
By Corollary~\ref{cor:strflat}, $\CC=A^e/M$ is a flat Fr\'echet $A^e$-module.
Applying Proposition~\ref{prop:flat-3} to the sequence $0\to A\to A_+\to\CC\to 0$
in $A^e\lmod(\Fr)$, we get the result.
\end{proof}

The following result generalizes Helemskii and Sheinberg's Theorem~\ref{thm:Hel_amen}
and gives a partial answer to a problem posed by Helemskii \cite[Problem 11]{X_31}.

\begin{theorem}
\label{thm:amen}
Let $A$ be a Fr\'echet-Arens-Michael algebra. The following conditions are equivalent:
\begin{enumerate}
\renewcommand{\theenumi}{\roman{enumi}}
\item $A$ is amenable;
\item for each Banach algebra $B$ such that there exists a homomorphism
$\varphi\colon A\to B$ with dense range, the algebra $B$ is amenable;
\item for each Arens-Michael decomposition
$A=\varprojlim (A_\lambda,\sigma^\mu_\lambda)_{\lambda\in\N}$, all the
$A_\lambda$'s are amenable Banach algebras;
\item $A$ is isomorphic to the reduced inverse limit of a sequence of
amenable Banach algebras;
\item $A_+$ is a strictly flat Fr\'echet $A$-bimodule;
\item $A$ is biflat and has a locally b.a.i.;
\item $I^\Delta$ has a right locally b.a.i.;
\item $A$ has a locally b.a.i., and $I^\Delta_0$ has a right locally b.a.i.
\end{enumerate}
If, in addition, $A$ is quasinormable, then the above conditions are equivalent
to the following:
\begin{enumerate}
\setcounter{enumi}{8}
\renewcommand{\theenumi}{\roman{enumi}}
\item $A$ is biflat and has a b.a.i.;
\item $I^\Delta$ has a right b.a.i.;
\item $A$ has a b.a.i., and $I^\Delta_0$ has a right b.a.i.
\end{enumerate}
\end{theorem}
\begin{proof}
$\mathrm{(i)}\Longrightarrow\mathrm{(ii)}$. This is a special case
of Proposition~\ref{prop:amen_dense}
(see also Remark~\ref{rem:amen_dense}).

$\mathrm{(ii)}\Longrightarrow\mathrm{(iii)}\Longrightarrow\mathrm{(iv)}$.
This is clear.

$\mathrm{(iv)}\Longrightarrow\mathrm{(vii)}$.
Let $A=\varprojlim (A_\lambda,\sigma^\mu_\lambda)_{\lambda\in\N}$ be a reduced inverse
limit of amenable Banach algebras. In what follows, we use the notation
introduced before Lemma~\ref{lemma:diagonal}.
By Helemskii and Sheinberg's Theorem~\ref{thm:Hel_amen},
each $I_\lambda^\Delta$ has a right b.a.i., which is equivalent to
(vii) by Lemma~\ref{lemma:diagonal} and Corollary~\ref{cor:bai_AMdec},

$\mathrm{(vii)}\iff\mathrm{(v)}\iff\mathrm{(i)}$.
By Lemma~\ref{lemma:diagonal},
for each $\lambda\in\N$ we have $I^\Delta_\lambda=\overline{\tau_\lambda(I^\Delta)}$,
so that $\overline{\tau_\lambda(I^\Delta)}$ is complemented in $A_\lambda^e$.
Now the equivalence of (vii), (v), and (i) follows from Corollary~\ref{cor:flat}
(cf. also Remark~\ref{rem:flat_unital}) applied to $A^e$ and $I^\Delta$.

$\mathrm{(iv)}\iff\mathrm{(viii)}$.
This is immediate from Helemskii and Sheinberg's Theorem~\ref{thm:Hel_amen},
Corollary~\ref{cor:bai_AMdec}, and Lemma~\ref{lemma:diagonal}.

$\mathrm{(i)}\&\mathrm{(viii)}\Longrightarrow\mathrm{(vi)}\Longrightarrow\mathrm{(i)}$.
This readily follows from Lemma~\ref{lemma:bifl_amen}.

Now assume that $A$ is quasinormable and satisfies conditions (i)--(viii).
Since quasinormability is inherited by projective tensor products
\cite[15.6.5]{Jarchow} and by quotients (and {\em a fortiori}
by complemented subspaces), we see that $A_0^e$, $A^e$, and
$I^\Delta$ are quasinormable. Together with
Theorem~\ref{thm:quasinorm}, this gives (ix) and (x).
Arguing as in \cite[Lemma VII.2.12 (III)]{X1}, we see that (ix) and (x)
together imply (xi).
\end{proof}

As an application of the above theorem, we can characterize
amenable Fr\'echet-Arens-Michael algebras in terms of derivations.

\begin{theorem}
Let $A$ be a Fr\'echet-Arens-Michael algebra. The following conditions are equivalent:
\begin{enumerate}
\renewcommand{\theenumi}{\roman{enumi}}
\item $A$ is amenable;
\item for each Banach $A$-bimodule $X$, every continuous derivation
from $A$ to $X^*$ is inner;
\item for each Fr\'echet $A$-bimodule $X$, every continuous derivation
from $A$ to $X^*$ is inner.
\end{enumerate}
\end{theorem}
\begin{proof}
$\mathrm{(i)}\Longrightarrow\mathrm{(iii)}$.
Let $X$ be a Fr\'echet $A$-bimodule.
Since $A_+$ is amenable, the augmented
complex~\eqref{std_hom_augm}
is exact. Therefore the dual complex of \eqref{std_hom_augm} is also exact
(see, e.g., \cite[26.4]{MV}). In particular, the short sequence
\begin{equation}
\label{std_cohom_012}
C_0(A,X)^* \to C_1(A,X)^* \to C_2(A,X)^*
\end{equation}
is exact. Using the same argument as in the Banach algebra case
(see \cite{X1}) and taking into account Proposition~\ref{prop:adj_ass},
we can identify \eqref{std_cohom_012} with the sequence
\begin{equation}
\label{H^1}
X^* \xra{\delta^0} \cL(A,X^*) \xra{\delta^1} \Bil(A\times A,X^*),
\end{equation}
where
\[
(\delta^0 f)(a)=a\cdot f-f\cdot a,\qquad
(\delta^1 \varphi)(a,b)=a\cdot\varphi(b)-\varphi(ab)+\varphi(a)\cdot b.
\]
Obviously, the kernel of $\delta^1$ is exactly the set of all continuous derivations
from $A$ to $X^*$, while the image of $\delta^0$ is exactly the set of all
inner derivations from $A$ to $X^*$. Since \eqref{H^1} is exact, this completes
the proof of (iii).

$\mathrm{(iii)}\Longrightarrow\mathrm{(ii)}$. This is clear.

$\mathrm{(ii)}\Longrightarrow\mathrm{(i)}$.
Let $A=\varprojlim A_\lambda$ be an Arens-Michael decomposition of $A$.
By Theorem~\ref{thm:amen}, it suffices to show that $A_\lambda$ is amenable for each
$\lambda$. Let $X$ be a Banach $A_\lambda$-bimodule, and let $D\colon A_\lambda\to X^*$
be a continuous derivation. Denoting by $\sigma_\lambda$ the canonical map
of $A$ to $A_\lambda$, we see that $D\sigma_\lambda\colon A\to X^*$ is a continuous
derivation. By (ii), $D\sigma_\lambda$ is inner. Since $\sigma_\lambda$ has
dense range, it follows that $D$ is also inner.
Thus $A_\lambda$ is amenable. In view of the above remarks, this
completes the proof of (i).
\end{proof}

Let us now turn to Johnson's Theorem~\ref{thm:Jhnsn}.

\begin{theorem}
\label{thm:amen_J}
Let $A$ be a Fr\'echet-Arens-Michael algebra. The following conditions are equivalent:
\begin{enumerate}
\renewcommand{\theenumi}{\roman{enumi}}
\item $A$ is amenable;
\item $A$ has a locally bounded approximate diagonal.
\end{enumerate}
If, in addition, $A$ is quasinormable, then the above conditions are equivalent
to the following:
\begin{enumerate}
\setcounter{enumi}{2}
\renewcommand{\theenumi}{\roman{enumi}}
\item $A$ has a bounded approximate diagonal;
\item $A$ has a virtual diagonal.
\end{enumerate}
\end{theorem}
\begin{proof}
$\mathrm{(i)}\iff\mathrm{(ii)}$. This is immediate from
Theorem~\ref{thm:amen} and Corollary~\ref{cor:bad_AMdec}.

Now assume that $A$ is amenable and quasinormable, so that it satisfies
condition (xi) of Theorem~\ref{thm:amen}.
Let $(e_\alpha)_{\alpha\in\Lambda}$
be a b.a.i. in $A$, and let $(u_\beta)_{\beta\in\Lambda'}$ be a right
b.a.i. in $I^\Delta_0$.
Define an order on $\Lambda\times\Lambda'$ by
$(\alpha,\beta)\prec (\gamma,\delta)$ iff $\alpha\prec\gamma$ and $\beta\prec\delta$.
For each $(\alpha,\beta)\in\Lambda\times\Lambda'$,
set $M_{\alpha\beta}=e_\alpha\otimes e_\alpha-u_\beta$.
We claim that $(M_{\alpha\beta})$ is a bounded approximate diagonal for $A$.
Indeed, it is clear that $\{ M_{\alpha\beta}\}$ is bounded.
Next, for each $a\in A$ we have
\begin{equation}
\label{bai_bad}
[a,M_{\alpha\beta}]=v(e_\alpha\otimes e_\alpha-u_\beta),
\end{equation}
where $v=a\otimes 1-1\otimes a\in I_0^\Delta$,
the product being taken in $A_0^e$. We have $vu_\beta\to v$,
as $(u_\beta)$ is a right a.i. in $I_0^\Delta$.
On the other hand, an easy computation (cf. \cite[VII.2.12]{X1},
\cite[2.9.21]{Dales}, \cite[8.1]{DW}, \cite[5.1.4]{Palmer} for the Banach algebra case)
shows that $(e_\alpha\otimes e_\alpha)$ is a b.a.i. in $A_0^e$.
Together with \eqref{bai_bad}, this implies that $[a,M_{\alpha\beta}]\to 0$
for each $a\in A$.

Furthermore, we have $\pi_0(M_{\alpha\beta})=e_\alpha^2$, and it is easily
seen that $(e_\alpha^2)$ is a b.a.i. in $A$. Indeed, for each continuous
submultiplicative seminorm $\|\cdot\|$ on $A$ we have
\[
\| e_\alpha^2 a-a\|
\le \| e_\alpha^2 a-e_\alpha a\| + \| e_\alpha a-a\|
\le \| e_\alpha a-a\| (\| e_\alpha\| + 1)\to 0,
\]
as $(e_\alpha)$ is bounded. Thus $\pi_0(M_{\alpha\beta})a\to a$ for each $a\in A$,
and so $(M_{\alpha\beta})$ is a bounded approximate diagonal for $A$.

$\mathrm{(iii)}\iff\mathrm{(iv)}$. This follows from Lemma~\ref{lemma:bad_vd}.
\end{proof}

Recall that a complete topological algebra $A$ with involution is a
{\em locally $C^*$-algebra} if the topology of $A$ can be defined by
a family of $C^*$-seminorms. A {\em $\sigma$-$C^*$-algebra}
(or a {\em Fr\'echet $C^*$-algebra}) is a metrizable locally $C^*$-algebra.
Each locally $C^*$-algebra is an Arens-Michael algebra. Moreover,
for each continuous $C^*$-seminorm
$p$ on $A$, the Banach algebra $A_p$ is a $C^*$-algebra in a natural way,
so that $A$ is a reduced inverse limit of $C^*$-algebras.
For details, see \cite{Frag_book}.

A locally $C^*$-algebra $A$ is {\em nuclear} \cite{Phil_JOT} if for each
continuous $C^*$-seminorm $p$ on $A$ the $C^*$-algebra $A_p$ is nuclear.
Equivalently, $A$ is nuclear if and only if $A$ is a reduced inverse
limit of nuclear $C^*$-algebras \cite[4.1]{Bh_Kar_nucl}.

Combining equivalence $\mathrm{(i)}\iff\mathrm{(ii)}$ from Theorem~\ref{thm:amen}
with the Connes-Haagerup Theorem (see \cite[6.5.12]{Runde}), we get the following.

\begin{corollary}
\label{cor:sigmacstar}
A $\sigma$-$C^*$-algebra is amenable if and only if
it is nuclear.
\end{corollary}

Following \cite{Mall_book}, we say that a Hausdorff topological
algebra $A$ is {\em uniform} if the topology on $A$ can be defined
by a directed family $\{\|\cdot\|_\lambda : \lambda\in\Lambda\}$ of
submultiplicative seminorms satisfying $\| a^2\|_\lambda=\| a\|^2_\lambda$
for each $a\in A$. Each uniform algebra is commutative and semisimple
\cite{Mall_book}.

Recall that a Hausdorff topological space $X$ is {\em hemicompact}
if the family of compact subsets of $X$ has a countable cofinal subset.
A Hausdorff topological space $X$ is a {\em $k$-space} if a
subset $F\subset X$ is closed whenever $F\cap K$ is closed for every
compact subset $K\subset X$.

Given a commutative locally convex algebra $A$, we denote by $\Omega(A)$
the spectrum of $A$ (i.e., the space of nonzero continuous characters)
endowed with the Gel'fand topology. In what follows, the algebra
$C(X)$ of continuous functions
on a topological space $X$ will be endowed with the compact-open topology.

\begin{corollary}
\label{cor:uniform}
Let $A$ be a unital uniform Fr\'echet algebra. The following conditions are
equivalent:
\begin{enumerate}
\renewcommand{\theenumi}{\roman{enumi}}
\item $A$ is amenable;
\item $A$ is topologically isomorphic to $C(X)$ for a
hemicompact $k$-space $X$;
\item the Gel'fand transform $\Gamma_A\colon A\to C(\Omega(A))$
is a topological isomorphism.
\end{enumerate}
\end{corollary}
\begin{proof}
$\mathrm{(i)}\Longrightarrow\mathrm{(iii)}$.
By \cite[4.1.3]{Gold_book}, $\Gamma_A$ is a topological isomorphism onto
its image, so we may identify $A$ with $\Gamma_A(A)$.
By \cite[3.2.8]{Gold_book}, $X=\Omega(A)$ is hemicompact, so that there exists
a countable exhaustion $X=\bigcup_n K_n$ with $K_n$ compact such that
each compact subset of $X$ is contained in some $K_n$.
Without loss of generality, we assume that $K_n\subset K_{n+1}$ for every $n$.
Let $A_n$ denote the uniform closure of the subalgebra
$\{ f|_{K_n} : f\in A\}$ of $C(K_n)$.
By \cite[4.1.6]{Gold_book}, the canonical map $A\to\varprojlim A_n$ is a topological
isomorphism. Using Theorem~\ref{thm:amen}, we see that $A_n$ is amenable
for each $n$. By Sheinberg's Theorem~\cite{Shnbrg_unif} (see also \cite[5.6.2]{Dales}),
$A_n=C(K_n)$. Thus we see that the composition of the maps
\[
A \hookrightarrow C(X) \to \varprojlim C(K_n)
\]
is an isomorphism. Since the second map is obviously injective, we conclude
that $A=C(X)$, as required.

$\mathrm{(iii)}\Longrightarrow\mathrm{(ii)}$.
This follows from \cite[3.1.9 (iii), Theorem]{Gold_book}.

$\mathrm{(ii)}\Longrightarrow\mathrm{(i)}$.
Let $X=\bigcup_n K_n$ be a countable exhaustion of $X$ with compact sets
having the above property. Then $C(X)\cong\varprojlim C(K_n)$ is an
Arens-Michael decomposition of $C(X)$ \cite[4.1.7]{Gold_book}.
Since each $C(K_n)$ is amenable (see, e.g., \cite[VII.2.40]{X1}, \cite[5.6.2]{Dales}),
Theorem~\ref{thm:amen} implies that $C(X)$ is amenable.
\end{proof}

We end this section by giving two examples of amenable Fr\'echet-Arens-Michael
algebras that are not covered by Corollaries~\ref{cor:sigmacstar}
and~\ref{cor:uniform}.

\begin{example}
Given a group $G$ and $g\in G$, let $\delta_g$ denote
the function on $G$ which is $1$ at $g$, $0$ elsewhere.
It is easy to see that for each group homomorphism $\varphi\colon G\to H$
we have a Banach algebra morphism $\varphi_*\colon\ell^1(G)\to\ell^1(H)$
uniquely determined by $\varphi_*(\delta_g)=\delta_{\varphi(g)}\; (g\in G)$.
Moreover, if $\varphi$ is onto, then so is $\varphi_*$.
Now let $\cG=(G_n,\varphi^m_n)_{n\in\N}$ be an inverse sequence of groups
such that the linking maps $\varphi^m_n$ are onto.
Set $\cL^1(\cG)=\varprojlim(\ell^1(G_n),(\varphi^m_n)_*)$.
Clearly, $\cL^1(\cG)$ is a Fr\'echet-Arens-Michael algebra.
A standard argument shows that $\cL^1(\cG)$ is not normable unless
the sequence $\cG$ stabilizes. Applying now Theorem~\ref{thm:amen},
we see that $\cL^1(\cG)$ is amenable if and only if all the $G_n$'s are
amenable.
\end{example}

Recall that a Fr\'echet space $E$ is a {\em quojection} if $E$
is isomorphic to a reduced inverse limit $\varprojlim(E_n,\varphi^m_n)$
of Banach spaces such that the linking maps $\varphi^m_n$ are onto.
All $\sigma$-$C^*$-algebras \cite[10.24]{Frag_book}
and all algebras of the form $\cL^1(\cG)$
are quojections. We now give an example of an amenable Fr\'echet-Arens-Michael
algebra which is not a quojection.

\begin{example}
Let $G$ be an infinite locally compact group, and let $A_p(G)\; (1<p<\infty)$
denote the Fig\`a-Talamanca-Herz algebra on $G$ (see, e.g., \cite[4.5]{Dales}).
The norm on $A_p(G)$ will be denoted by $\|\cdot\|_p$.
Suppose that $G$ is amenable. Then, by \cite{Herz1}, for each $p\ge q\ge 2$
we have $A_q(G)\subset A_p(G)$, and the inclusion of $A_q(G)$
into $A_p(G)$ is continuous with dense range. Moreover,
as was shown in \cite{Cowl_Fourn}, $A_p(G)\ne A_q(G)$ if $p\ne q\ge 2$.
Given $p\ge 2$, set $A_{p+}(G)=\bigcap_{q>p} A_q(G)=\varprojlim_{q>p} A_q(G)$.
Clearly, $A_{p+}(G)$ is a Fr\'echet-Arens-Michael algebra.
To show that $A_{p+}(G)$ is not a quojection, fix any $q>p$ and note that
$\|\cdot\|_q$ is a continuous norm on $A_{p+}(G)$. By \cite[8.4.33]{PCB},
a quojection having a continuous norm is normable.
Therefore if $A_{p+}(G)$ were a quojection, the norms $\|\cdot\|_q$ and
$\|\cdot\|_r$ on $A_{p+}(G)$ would be equivalent whenever $q,r$ are close enough to $p$,
which contradicts the fact that the inclusion $A_q(G)\subset A_r(G)\; (2\le q<r)$
is proper.

Now assume that $G$ contains an abelian subgroup of finite index.
Then $A_2(G)$ is amenable \cite{LLW}, and hence $A_p(G)$ is amenable
for every $p$, since the inclusion $A_2(G)\to A_p(G)$ has dense range.
Now Theorem~\ref{thm:amen} implies that $A_{p+}(G)$ is amenable.
\end{example}

\section{A counterexample}
\label{sect:counterexample}
Now we present an example of a commutative
Fr\'echet-Arens-Michael algebra with a locally b.a.i., but
without a b.a.i. Together with Corollary~\ref{cor:strflat}, this will show
that Helemskii and Sheinberg's Theorem~\ref{thm:Hel_flat}
does not extend {\em verbatim} to Fr\'echet-Arens-Michael
algebras.

Let $P$ be a family of real-valued sequences such that for each $i\in\N$
there exists $p\in P$ with $p_i>0$. Suppose also that $P$ is directed, i.e.,
for each $p,q\in P$ there exists $r\in P$ such that $r_i\ge\max\{ p_i,q_i\}$
for all $i\in\N$. Recall that the {\em K\"othe sequence space} $\lambda(P)$
is defined as follows:
\[
\lambda(P)=\Bigl\{ a=(a_i)_{i\in\N}\in\CC^\N :
\| a\|_p=\sum_i |a_i| p_i < \infty\;\forall p\in P\Bigr\}.
\]
In the sequel, for each $i\in\N$ we set $e_i=(0,\ldots,0,1,0,\ldots)$,
where the single nonzero entry is in the $i$th slot.
The linear span of the $e_i$'s is denoted by $c_{00}$.

\begin{lemma}
\label{lemma:A(P)}
Suppose that $p_i\ge 1$ for each $i\in\N$. Then
there exists a unique continuous product on $\lambda(P)$
such that $e_i e_j=e_{\min\{ i,j\}}$ for all $i,j\in\N$.
Together with this product, $\lambda(P)$ becomes a commutative Arens-Michael algebra.
\end{lemma}
\begin{proof}
Clearly, the map $(e_i,e_j)\mapsto e_{\min\{ i,j\}}$ uniquely extends to an associative
product on $c_{00}$.
Since $p_i\ge 1$, we have $\| e_i e_j\|_p\le \| e_i \|_p \| e_j\|_p$ for all $i,j\in\N$.
This easily implies that $\| ab\|_p \le \| a\|_p \| b\|_p$ for all $a,b\in c_{00}$.
Using the density of $c_{00}$ in $\lambda(P)$ and extending the product by continuity
to the whole of $\lambda(P)$, we get the result.
\end{proof}

The Arens-Michael algebra defined in Lemma~\ref{lemma:A(P)} will be denoted by $A(P)$.
Note that for each $p\in P$ the Banach algebra $A(P)_p$ is
isometrically isomorphic to $A(\{ p\})$.

\begin{lemma}
\label{lemma:A(P)_nonunit}
The algebra $A(P)$ has no identity.
\end{lemma}
\begin{proof}
Assume, towards a contradiction, that $u=(u_i)$ is an identity in $A(P)$.
Then for each $j\in\N$ we have
\[
e_j=e_j u = \sum_i u_i e_j e_i=\sum_{i<j} u_i e_i+\Bigl(\sum_{i\ge j} u_i\Bigr)e_j.
\]
This implies that $u_i=0$ for all $i<j$. As $j\in\N$ was arbitrary, we obtain $u=0$,
which is a contradiction.
\end{proof}

\begin{lemma}
\label{lemma:A(P)_lbai}
Suppose that each sequence $p\in P$ has a bounded subsequence. Then $A(P)$ has a
locally b.a.i.
\end{lemma}
\begin{proof}
Let $p\in P$, and let $(p_{n_k})_{k\in\N}$ be a bounded subsequence of $p$.
Clearly, for each $a\in c_{00}$ we have $a e_{n_k}=a$ for $k$ large enough.
Since $\| e_{n_k}\|=p_{n_k}$, this implies that $e_{n_k}$ is a b.a.i.
in $(c_{00},\| \cdot\|_p)$ and hence in $A(\{ p\})$ (cf. \cite[1.4]{DW}).
Now the result follows from Corollary~\ref{cor:bai_AMdec}.
\end{proof}

Now we define $P$ exactly as in the classical K\"othe-Grothendieck example
\cite{Kothe_stufen,Groth_F_DF}.
Namely, for each $k\in\N$ we define an infinite matrix
$\alpha^{(k)}=(\alpha^{(k)}_{ij})_{i,j\in\N}$
by setting
$$
\alpha^{(k)}_{ij}=\left\{
\begin{array}{ll}
j^k, & i<k\\
i^k, & i\ge k.
\end{array}
\right.
$$
Fix a bijection $\varphi\colon\N^2\to\N$ such that $\varphi(i,j+1)<\varphi(i,j)$
for all $i,j\in\N$. For each $k\in\N$ define a sequence $p^{(k)}=(p^{(k)}_n)_{n\in\N}$
by $p^{(k)}_n=\alpha^{(k)}_{\varphi^{-1}(n)}$. Finally, set $P=\{ p^{(k)} : k\in\N\}$.
Since $P$ is countable, $A(P)$ is a Fr\'echet algebra.

\begin{prop}
If $P$ is as above, then $A(P)$ has a locally b.a.i., but does not have a b.a.i.
\end{prop}
\begin{proof}
Given $p=p^{(k)}\in P$, set $n_i=\varphi(k,i)\; (i\in\N)$.
By assumption, the sequence $(n_i)$ is strictly increasing.
We have $p_{n_i}^{(k)}=\alpha_{ki}^{(k)}=k^k$, so that $(p_{n_i}^{(k)})_{i\in\N}$
is a bounded subsequence of $p^{(k)}$. By Lemma~\ref{lemma:A(P)_lbai},
$A(P)$ has a locally b.a.i. On the other hand, $A(P)$ is a Montel
space \cite{Kothe_stufen,Groth_F_DF}. Applying Lemma~\ref{lemma:bai_Montel}
and Lemma~\ref{lemma:A(P)_nonunit}, we see that $A(P)$ does not have a b.a.i.
\end{proof}

Together with Corollary~\ref{cor:strflat} this gives the following.

\begin{corollary}
There exists a commutative Fr\'echet-Arens-Michael algebra $A$
such that the trivial Fr\'echet $A$-module $\CC=A_+/A$ is strictly flat,
but $A$ does not have a b.a.i. In addition, the underlying Fr\'echet
space of $A$ is Montel.
\end{corollary}

\noindent\textbf{Acknowledgments. }The author is grateful to
A.~Ya.~Helemskii for helpful discussions.


\begin{thebibliography}{99}
\bibitem{Allan_stabseq}
Allan, G. R.
{\em Stable inverse-limit sequences, with applications to Fr\'echet algebras}.
Studia Math. \textbf{121} (1996), 277--308.
\bibitem{Allan_auto}
Allan, G. R.
{\em Stable inverse-limit sequences and automatic continuity}.
Studia Math. \textbf{141} (2000), 99--107.
\bibitem{Arens_gennorm}
Arens, R.
{\em A generalization of normed rings}.
Pacific J. Math. \textbf{2} (1952), 455--471.
\bibitem{Bh_Kar_nucl}
Bhatt, S. J.; Karia, D. J.
{\em Complete positivity, tensor products and $C^*$-nuclearity for
inverse limits of $C\sp *$-algebras}.
Proc. Indian Acad. Sci. Math. Sci. \textbf{101} (1991), no.~3, 149--167.
\bibitem{Cowl_Fourn}
Cowling, M. G.; Fournier, J. J. F.
{\em Inclusion and noninclusion of spaces of convolution operators},
Trans. Amer. Math. Soc. \textbf{221} (1976), no.~1, 59--95.
\bibitem{Dales}
Dales, H. G.
{\em Banach Algebras and Automatic Continuity},
Clarendon Press, Oxford, 2000.
\bibitem{DLZ}
Dales, H. G.; Loy, R. J.; Zhang, Y.
{\em Approximate amenability for Banach sequence algebras}.
Studia Math., to appear.
\bibitem{DW}
Doran, R. S.; Wichmann, J.
{\em Approximate Identities and Factorization in Banach Modules},
Lecture Notes in Math. \textbf{768},
Springer-Verlag, Berlin--New York, 1979.
\bibitem{Eng}
Engelking, R.
{\em General Topology},
PWN, Warsaw, 1977.
\bibitem{EschmPut}
Eschmeier, J.; Putinar, M.
{\itshape Spectral Decompositions and Analytic Sheaves}.
Clarendon Press, Oxford, 1996.
\bibitem{Frag_book}
Fragoulopoulou, M.
{\em Topological Algebras with Involution},
North-Holland Mathematics Studies, 200, Elsevier, Amsterdam, 2005.
\bibitem{GL}
Ghahramani, F.; Loy, R. J.
{\em Generalized notions of amenability}.
J. Funct. Anal. \textbf{208} (2004), no.~1, 229--260.
\bibitem{Gold_book}
Goldmann, H.
{\em Uniform Fr\'echet Algebras}.
North-Holland Mathematics Studies, 162.
North-Holland, Amsterdam--New York, 1990.
\bibitem{Groth_F_DF}
Grothendieck, A.
{\em Sur les espaces (F) et (DF)},
Summa Brasil. Math. \textbf{3} (1954), 57--123.
\bibitem{Guichardet1}
Guichardet, A.
{\itshape Sur l'homologie et la cohomologie des alg\`ebres de Banach},
C. R. Acad. Sci. Paris Ser. A \textbf{262} (1966), 38--41.
\bibitem{Gulick}
Gulick, S. L.
{\em The bidual of a locally multiplicatively-convex algebra},
Pacific J. Math. \textbf{17} (1966), 71--96.
\bibitem{Hel_period}
Helemskii, A. Ya.
{\em A periodic product of modules over Banach algebras} (Russian),
Funkcional. Anal. i Prilo\v zen \textbf{5} (1971), no.~1, 95--96.
English transl.: Functional Anal. Appl. \textbf{5} (1971), 84--85.
\bibitem{Hel_Shein}
Helemskii, A. Ya.; Sheinberg, M. V.
{\em Amenable Banach algebras} (Russian).
Funktsional. Anal. i Prilozhen. \textbf{13} (1979), no.~1, 42--48.
English transl.: Functional Anal. Appl. \textbf{13} (1979), no.~1, 32--37.
\bibitem{Hel_flat_amen}
Helemskii, A. Ya.
{\em Flat Banach modules and amenable algebras} (Russian).
Trudy Moskov. Mat. Obshch. \textbf{47} (1984), 179--218.
English transl.: Trans. Moscow Math. Soc. (1985), 199-244.
\bibitem{X1}
Helemskii, A. Ya. {\itshape The Homology of Banach and Topological Algebras},
Moscow University Press, 1986 (Russian); English transl.: Kluwer Academic
Publishers, Dordrecht, 1989.
\bibitem{X2}
Helemskii, A. Ya. {\itshape Banach and Polynormed Algebras: General Theory,
Representations, Homology}, Nauka, Moscow, 1989 (Russian); English transl.:
Oxford University Press, 1993.
\bibitem{X_31}
Helemskii, A. Ya.
{\itshape 31 problems of the homology of the algebras of analysis},
``Linear and Complex Analysis: Problem Book~3, Part~I''
(eds. V.~P.~Havin and N.~K.~Nikolski),
Lecture Notes in Math. \textbf{1573} (1994), 54--78, Berlin, Springer.
\bibitem{X_HOA}
Helemskii, A. Ya.
{\itshape Homology for the Algebras of Analysis},
Handbook of algebra, Vol.~2 (ed. M.~Hazewinkel), 151--274,
Amsterdam, North-Holland, 2000.
\bibitem{Herz1}
Herz, C.
{\em The theory of $p$-spaces with an application to convolution operators},
Trans. Amer. Math. Soc \textbf{154} (1971), 69--82.
\bibitem{Jarchow}
Jarchow, H.
{\em Locally Convex Spaces}.
Teubner, Stuttgart, 1981.
\bibitem{Jhnsn_centr}
Johnson, B. E.
{\em An introduction to the theory of centralizers},
Proc. London Math. Soc. (3) \textbf{14} (1964), 299--320.
\bibitem{Jhnsn_CBA}
Johnson, B. E.
{\em Cohomology in Banach algebras},
Mem. Amer. Math. Soc. \textbf{127} (1972).
\bibitem{Jhnsn_appr}
Johnson, B. E.
{\em Approximate diagonals and cohomology of certain annihilator
Banach algebras}.
Amer. J. Math. \textbf{94} (1972). no.~3, 685--698.
\bibitem{Kothe_stufen}
K\"othe, G.
{\em Die Stufenr\"aume, eine einfache Klasse linearer vollkommener R\"aume}.
Math. Z. \textbf{51} (1948), 317--345.
\bibitem{Kothe2}
K\"othe, G.
{\itshape Topological Vector Spaces} II,
Springer-Verlag, New York, 1979.
\bibitem{LLW}
Lau, A. T.-M.; Loy, R. J.; Willis, G. A.
{\em Amenability of Banach and $C^*$-algebras on locally compact groups},
Studia Math. \textbf{119} (1996), no.~2, 161--178.
\bibitem{MacLane}
MacLane, S.
{\em Homology}.
Springer-Verlag, Berlin, 1963.
\bibitem{ML_work}
MacLane, S.
{\em Categories for the Working Mathematician},
Graduate Texts in Mathematics, Vol. 5.
Springer-Verlag, New York-Berlin, 1971.
\bibitem{Mall_book}
Mallios, A.
{\em Topological Algebras. Selected Topics}.
North-Holland Mathematics Studies, 124.
North-Holland, Amsterdam--New York, 1986.
\bibitem{MV}
Meise, R.; Vogt, D.
{\em Introduction to Functional Analysis},
Clarendon Press, Oxford, 1997.
\bibitem{Michael}
Michael, E. A.
{\em Locally multiplicatively-convex topological algebras},
Mem. Amer. Math. Soc. \textbf{11} (1952).
\bibitem{Pal_projlim}
Palamodov, V. P.
{\em The projective limit functor in the category of topological linear spaces} (Russian).
Mat. Sb. (N.S.) \textbf{75} (1968), 567--603.
English transl.: Math. USSR-Sb. \textbf{75} (1968), 529--558.
\bibitem{Pal_methods}
Palamodov, V. P.
{\em Homological methods in the theory of locally convex spaces} (Russian),
Uspehi Mat. Nauk \textbf{26} (1971), no.~1 (157), 3--65.
English transl.: Russian Math. Surveys \textbf{26} (1971), 1--64.
\bibitem{Palmer}
Palmer, T. W.
{\em Banach Algebras and the General Theory of $*$-Algebras. Vol. I.
Algebras and Banach algebras}, Encyclopedia of Mathematics and its
Applications, 49, Cambridge University Press, Cambridge, 1994.
\bibitem{PCB}
P\'erez Carreras, P.; Bonet, J.
{\em Barrelled Locally Convex Spaces},
North-Holland Mathematics Studies, 131.
North-Holland, Amsterdam--New York, 1987.
\bibitem{Phil_JOT}
Phillips, N. C.
{\em Inverse limits of $C\sp *$-algebras}.
J. Operator Theory \textbf{19} (1988), no.~1, 159--195.
\bibitem{Pir_inject}
Pirkovskii, A. Yu.
{\itshape On the problem of the existence of a sufficient
number of injective Fr\'echet modules over nonnormable
Fr\'echet algebras} (in Russian),
Izvestiya Ross. Akad. Nauk, ser. matem.
\textbf{62} (1998), no.~4, 137--154. English transl.:
Izvestiya: Mathematics \textbf{62} (1998), no.~4, 773--788.
\bibitem{Pir_injdim}
Pirkovskii, A. Yu.
{\itshape On Arens-Michael algebras which
do not have non-zero injective $\widehat\otimes$-modules},
Studia Math. \textbf{133} (1999), no.~2, 163--174.
\bibitem{Pir_nova}
Pirkovskii, A. Yu.
{\em Injective topological modules, additivity formulas for homological
dimensions and related topics},
Topological Homology: Helemskii's Moscow Seminar (ed. A.~Ya.~Helemskii), 93--143,
Nova Sci. Publ., Huntington, NY, 2000.
\bibitem{Pir_ICTAA05}
Pirkovskii, A. Yu.
{\em Strictly flat cyclic Fr\'echet modules and approximate identities}.
Proc. of the 5th International Conference on Topological Algebras and
Applications, Athens, 2005. Contemporary Math., to appear.
Preprint arXiv.org:math.FA/0511132.
\bibitem{Pir_qfree}
Pirkovskii, A. Yu.
{\em Arens-Michael envelopes, homological epimorphisms, and
relatively quasi-free algebras}. Trans. Moscow Math. Soc., to appear.
\bibitem{Podara}
Podara, C. P.
{\em On strictly flat Fr\'echet modules}.
Proc. of the 5th International Conference on Topological Algebras and
Applications, Athens, 2005. Contemporary Math., to appear.
\bibitem{Prosm_DFA}
Prosmans, F.
{\em Derived categories for Functional Analysis}.
Publ. Res. Inst. Math. Sci. \textbf{36} (2000), no.~1, 19--83.
\bibitem{Quillen}
Quillen, D.
{\em Higher algebraic $K$-theory I}.
Lecture Notes in Math. \textbf{341}, 85--147,
Berlin, Springer, 1973.
\bibitem{Runde}
Runde, V.
{\em Lectures on Amenability},
Lecture Notes in Math. 1774, Springer, 2002.
\bibitem{Sch}
Schaefer, H. {\itshape Topological Vector Spaces},
Macmillan, New York, 1966.
\bibitem{Sel_cohchar}
Selivanov, Yu. V.
{\em Cohomological characterizations of biprojective and biflat Banach algebras}.
Monatsh. Math. \textbf{128} (1999), no.~1, 35--60.
\bibitem{Shnbrg_unif}
Sheinberg, M. V.
{\em A characterization of the algebra $C(\Omega)$ in terms of cohomology groups}
(Russian). Uspekhi Mat. Nauk \textbf{32} (1977), 203--204.
\bibitem{Stnstrm}
Stenstr{\"o}m, B.
{\em Rings of Quotients}.
Die Grundlehren der Mathematischen Wissenschaften, Band 217.
Springer-Verlag, New York, 1975.
\bibitem{T1}
Taylor, J. L. {\itshape Homology and cohomology for topological algebras},
Adv. Math. \textbf{9} (1972), 137--182.
\bibitem{Wngnrth}
Wengenroth, J.
{\em Derived Functors in Functional Analysis},
Lecture Notes in Math. \textbf{1810}, Berlin, Springer, 2003.

\end{thebibliography}
\end{document}